\documentclass[a4paper]{article}
\usepackage{amsmath,amsthm,amssymb,bm,Fermionic_form}
\usepackage[all]{xy}

\begin{document}
\title{Fusion products of Kirillov-Reshetikhin modules and the $X=M$ conjecture}
\author{Katsuyuki Naoi}
\date{}

\maketitle

\begin{abstract}
  In this article, we show in the $ADE$ case that the fusion product of Kirillov-Reshetikhin modules for a current algebra,
  whose character is expressed in terms of fermionic forms, 
  can be constructed from one-dimensional modules by using Joseph functors.
  As a consequence, we obtain some identity between fermionic forms and Demazure operators.
  Since the same identity is also known to hold for one-dimensional sums of nonexceptional type,
  we can show from these results the $X=M$ conjecture for type $A_n^{(1)}$ and $D_n^{(1)}$.
\end{abstract}

\section{Introduction}

Let $\fg$ be an affine Kac-Moody Lie algebra with index set $I$, 
and $U_q'(\fg)$ the corresponding quantum affine algebra without the degree
operator.
In \cite{MR1745263,MR1903978}, it was conjectured that a certain subfamily of finite-dimensional $U_q'(\fg)$-modules
known as Kirillov-Reshetikhin (KR for short) modules $W^{r,\ell}$ has crystal bases $B^{r,\ell}$ called 
KR crystals.
Here the index $r$ corresponds to a node of $I_0 = I \setminus \{ 0\}$ where $0$ is the node specified in \cite{MR1104219}, 
and $\ell$ is a positive integer.
This conjecture has been confirmed in many cases, in particular when $\fg$ is nonexceptional \cite{MR2403558}.

Let $B = B^{r_p,\ell_p} \otimes \cdots \otimes B^{r_1,\ell_1}$ be a tensor product of KR crystals.
In \cite{MR1745263,MR1903978}, the authors defined the one-dimensional sum
\[ X(B, \mu ,q) = \sum_{\begin{smallmatrix} b \in B \\ \wti{e}_i b= 0 \ (i \in I_0) \\ \wt(b) = 
   \mu \end{smallmatrix}} q^{D(b)}\in \Z[q,q^{-1}],
\]
where $\mu$ is a dominant integral weight of the simple Lie subalgebra $\fg_0 \subseteq \fg$ 
whose Dynkin nodes are $I_0$,
$D$ is a certain $\Z$-function on $B$ called the energy function, and $\wti{e}_i$ are Kashiwara operators.
Then they conjectured that $X(B,\mu,q)$ has an explicit expression $M(\bm{\nu},\mu,q)$ called the fermionic form,
where $\bm{\nu}$ is the sequence of elements of $I_0 \times \Z_{> 0}$ corresponding to $B$.
This conjecture is called the $X=M$ conjecture,
which was confirmed in many instances \cite{MR1890195,MR2119115,MR2194969}, 
but a proof in full generality has not been available 
except for the type $\fg=A_n^{(1)}$ in \cite{MR1890195}.
It should be mentioned that this conjecture was recently settled for nonexceptional $\fg$ 
under the large rank hypothesis by combining the results in \cite{MR1890195}, \cite{LOS} and \cite{OS}, 
and for general $\fg$ if $\ell_1 = \cdots = \ell_p =1$ \cite{Na}.

Recently, an identity was proved in \cite{Na3} which connects one-dimensional sums 
of nonexceptional type with Demazure operators.
Let us recall the result briefly.
For simplicity, we assume $\fg$ is of type $A_n^{(1)}$ or $D_n^{(1)}$ here.
Denote by $P$ the weight lattice of $\fg$
and by $D_i$ for $i \in I$ the Demazure operator on the group ring $\Z[P]$ defined by
\begin{equation*}
   D_i(f) = \frac{f - e^{-\ga_i}\cdot s_i(f)}{1-e^{-\ga_i}},
\end{equation*}
where the simple reflection $s_i$ acts on $\Z[P]$ by $s_i(e^{\gl}) = e^{s_i(\gl)}$.
Let $W$ be the Weyl group of $\fg$, $\wti{W}$ the extended affine Weyl group, 
and $\gS$ the subgroup of the group of Dynkin automorphisms 
such that $\wti{W} \cong W\rtimes \gS$ (for the precise definition, see Section \ref{Nontwisted affine Lie algebra}).
Then for each $\tau \in \gS$ and $w \in W$ with a reduced expression $w = s_{i_k} \cdots s_{i_1}$,
$D_{w \tau}$ is defined by $D_{w\tau} = D_{i_k} \cdots D_{i_1} \circ \tau$
(this definition does not depend on the choice of the expression). 
Let now $B= B^{r_p,\ell_p} \otimes \cdots \otimes B^{r_1,\ell_1}$
be a tensor product of KR crystals such that $\ell_1 \le \cdots \le \ell_p$,
$P_0^+$ the set of dominant integral weights of $\fg_0$, and $V_{\fg_0}(\mu)$ the irreducible 
$\fg_0$-module with highest weight $\mu \in P_0^+$.
Then the following identity was proved in \cite{Na3}, where $C$ is some constant and we set $q = e^{-\gd}$:
\begin{align}\label{eq:intro}
  q^{C}&e^{\ell_p\gL_0}\sum_{\mu \in P_0^+} X(B, \mu,q) \, \ch V_{\fg_0}(\mu) \\
  & =D_{t_{w_0(\varpi_{r_p})}}\Big( e^{(\ell_p-\ell_{p-1})\gL_0}
  \cdots  D_{t_{w_0(\varpi_{r_2})}}
  \big( e^{(\ell_2 -\ell_1)\gL_0}\cdot D_{t_{w_0(\varpi_{r_1})}}(e^{\ell_1\gL_0})\big)\cdots \Big). \nonumber
\end{align}
Here $\gd$ is the null root, $\gL_0$ is the fundamental weight of $\fg$ associated with the node $0$,
$\varpi_i$ are the fundamental weights of $\fg_0$, $w_0$ is the longest element of the Weyl group of $\fg_0$,
and $t_{w_0(\varpi_{r_j})} \in \wti{W}$ are the translations.

The goal of this article is to show the identity (\ref{eq:intro}) with the one-dimensional sums replaced 
by the corresponding fermionic forms.
Namely, we will prove the following result as a corollary of Theorem \ref{Thm:intro} stated below (Corollary \ref{Cor:M-side}):

\begin{Cor}\label{Cor:M-side_in_intro}
  Assume that $\fg$ is of nontwisted type and $\fg_0$ is of $ADE$ type. 
  Let $\bm{\nu}= \big((r_1,\ell_1),\ldots, (r_p,\ell_p)\big)$ be a sequence of elements of 
  $I_0 \times \Z_{>0}$ such that $\ell_1 \le \cdots \le \ell_p$. 
  Then we have
  \begin{align*}
    q^C&e^{\ell_p\gL_0}\sum_{\mu \in P_0^+} M(\bm{\nu}, \mu,q) \ch V_{\fg_0}(\mu)\\ 
    &= D_{t_{w_0(\varpi_{r_p})}}\Big(e^{(\ell_p-\ell_{p-1})\gL_0} \cdots
    D_{t_{w_0(\varpi_{r_2})}}
    \big( e^{(\ell_2-\ell_1)\gL_0} \cdot D_{t_{w_0(\varpi_{r_1})}}(e^{\ell_1\gL_0})\big) \cdots \Big)
  \end{align*}
  with some constant $C$, where we set $q = e^{-\gd}$.
\end{Cor}

Then as an immediate consequence of (\ref{eq:intro}) and this corollary, the $X=M$ conjecture for $A_n^{(1)}$ and $D_n^{(1)}$
is settled (Theorem \ref{Thm:X=M_conj}).

Our strategy of the proof of Corollary \ref{Cor:M-side_in_intro} is to show the isomorphism between two modules
whose characters are equal to the left hand side and the right hand side respectively.
Let us describe what these modules are.
The module corresponding to the left hand side is the fusion product of Kirillov-Reshetikhin (KR)
modules for the current algebra $\fg_0 \otimes \C[t]$.
KR modules for $\fg_0 \otimes \C[t]$, which we denote by $KR^{r,\ell}$ ($r \in I_0$, $\ell \in \Z_{>0}$) in this article,
were defined in \cite{MR1836791,MR2238884} in terms of generators and relations.
They can also be obtained from $W^{r,\ell}$ (KR modules for $U_q'(\fg)$) by specialization and restriction,
which is why they are so named.
The fusion product is a refinement of the usual tensor product defined in \cite{MR1729359}, 
which constructs a cyclic graded ($\fg_0 \otimes \C[t]$)-module from some such modules.
It was proved by Di Francesco and Kedem in \cite{MR2428305}
that the fusion product of KR modules has the required (graded) character.
Namely, they proved the following character formula for the fusion product of KR modules:
for a sequence $\bm{\nu} = \big((r_1,\ell_1), \ldots ,(r_p,\ell_p)\big)$ of elements of $I_0 \times \Z_{>0}$,
we have 
\begin{equation}\label{eq:identity2}
  \ch KR^{r_p,\ell_p} * \cdots * KR^{r_1, \ell_1} = \sum_{\mu \in P_0^+} M(\bm{\nu}, \mu,q) \ch V_{\fg_0}(\mu),
\end{equation}
where $*$ denotes the fusion product.

To present the module corresponding to the right hand side, we recall the definition of Joseph functors
introduced by Joseph \cite{MR826100}.
Let $\fh$ be a Cartan subalgebra of $\fg$, $\fb$ a Borel subalgebra of $\fg$ containing $\fh$, 
and $\fp_i$ for $i \in I$ the parabolic subalgebra associated with the set $\{i\}$ containing $\fb$.
We denote by $\Mod{\fb}$ (resp.\ $\Mod{\fp_i}$) the category of finite-dimensional $\fh$-semisimple left 
$\fb$-modules (resp.\ $\fp_i$-modules).
Then the functor $\cD_i\colon \Mod{\fb} \to \Mod{\fp_i}$ is defined as the left adjoint functor of the restriction functor
$\Mod{\fp_i} \to \Mod{\fb}$.
For $w \in W$ with a reduced expression $w= s_{i_k} \cdots s_{i_1}$, the Joseph functor $\cD_w$
is defined by 
\[ \cD_w= \cD_{i_k} \cdots \cD_{i_1}\colon \Mod{\fb} \to \Mod{\fb},
\]
which does not depend on the choice of the expression.
We also define $\cD_{w\tau}$ for $w \in W$ and $\tau \in \gS$ by $\cD_{w\tau}= \cD_w \circ (\tau^{-1})^*$, 
where $(\tau^{-1})^*$ denotes
the pull-back functor associated with the Lie algebra automorphism $\tau^{-1}$ of $\fb$.
Then the $\fb$-module corresponding to the right hand side of Corollary \ref{Cor:M-side_in_intro}
is constructed using these functors as follows:
\[ \cD_{t_{w_0(\varpi_{r_p})}}\Big(\C_{(\ell_p-\ell_{p-1})\gL_0} \otimes \cdots \otimes \cD_{t_{w_0(\varpi_{r_2})}}
     \big( \C_{(\ell_2-\ell_1)\gL_0} \otimes \cD_{t_{w_0(\varpi_{r_1})}}\C_{\ell_1\gL_0}\big) \cdots \Big).
\]
From the results of \cite{MR826100} and \cite{MR1923198},
we can easily see that $\ch \cD_w M = D_w\,\ch M$ holds for $w \in \wti{W}$ if a $\fb$-module $M$ has 
a Demazure flag (see Definition \ref{Def:Demazure_flag}).
It is checked from this fact that the above module in fact has the required character.

Now, what we have to prove is the following isomorphism, which is our main theorem in this article (Theorem 
\ref{Thm:Main_theorem}):

\begin{Thm}\label{Thm:intro}
  Assume that $\fg$ is of nontwisted type and $\fg_0$ is of type $ADE$.
  Let $KR^{r_1,\ell_1}, \ldots, KR^{r_p,\ell_p}$ be a sequence of KR modules such that $\ell_1 \le \cdots \le \ell_p$.
  Then there exists an isomorphism of $\fb$-modules
  \begin{align*}
     &\C_{\ell_p\gL_0 + C \gd}\otimes \left(KR^{r_p,\ell_p} * \cdots *KR^{r_2,\ell_2}*KR^{r_1,\ell_1}\right) \nonumber \\ 
     & \
     \cong\cD_{t_{w_0(\varpi_{r_p})}}\Big(\C_{(\ell_p-\ell_{p-1})\gL_0} \otimes \cdots \otimes \cD_{t_{w_0(\varpi_{r_2})}}
     \big( \C_{(\ell_2-\ell_1)\gL_0} \otimes \cD_{t_{w_0(\varpi_{r_1})}}\C_{\ell_1\gL_0}\big) \cdots \Big)
  \end{align*}
  with some constant $C$, where the left hand side is naturally regarded as a $\fb$-module.
\end{Thm}

As explained above, Corollary \ref{Cor:M-side_in_intro} follows as a direct consequence of this theorem.
It should be remarked that this theorem for $\fg = A_1^{(1)}$ is already proved by Feigin and Loktev in \cite{MR1729359}.

The plan of this article is as follows.
In Section \ref{Section:Demazure}, we review the results on Demazure modules and Joseph functors.
In Section \ref{Nontwisted affine Lie algebra}, we prepare some notation and elementary lemmas concerning a nontwisted affine 
Lie algebra.
In Section \ref{Section:KR}, we recall the definition of KR modules for a current algebra.
For the later convenience, we define them as $\fp_{I_0}$-modules in this article,
where $\fp_{I_0}$ denotes the parabolic subalgebra associated with $I_0$.
In Section \ref{Fusion product}, we recall the definition of fusion products.
Then we state our main theorem in Section \ref{Section:Main}, and explain in Section \ref{X=M_conjecture}
how the $X=M$ conjecture for type $A_n^{(1)}$ and $D_n^{(1)}$ is deduced from this theorem.
The last two sections are devoted to prove the main theorem.
In Section \ref{Section:b-fusion}, we introduce the $\fb$-fusion product, which is some modified version of the fusion product
constructing a $\fb$-module from some $\fb$-modules. 
Then in Section \ref{Section:proof}, we give the proof of the main theorem using $\fb$-fusion products.

\ \\
\textbf{Acknowledgments:} 
The author would like to express his gratitude to R.\ Kodera, S.\ Naito and Y.\ Saito
for a lot of helpful discussions and comments.

\section{Demazure modules and Joseph functors}\label{Section:Demazure}

Let $\fg$ be a complex symmetrizable Kac-Moody Lie algebra with index set $I$ and 
Chevalley generators $\{e_i, f_i \mid i \in I\}$, 
$\fb \subseteq \fg$ its Borel subalgebra, and $\fh \subseteq \fb$ its Cartan subalgebra.
Let $\ga_i \in \fh^*$ $(i \in I)$ be the simple roots, $\ga_i^{\vee} \in \fh$ $(i \in I)$ the simple coroots,
and $W$ the Weyl group of $\fg$ with simple reflections $s_i$ ($i \in I$).
Denote by 
\begin{align*} 
   P &= \{ \gl \in \fh^* \mid \langle \gl, \ga_i^{\vee} \rangle \in \Z \ \text{for $i \in I$} \} \ \text{and} \\ 
   P^+&= \{ \gl \in P \mid \langle \gl, \ga_i^{\vee} \rangle \in \Z_{\ge 0} \ \text{for $i \in I$} \}
\end{align*}
the weight lattice and the set of dominant integral weights respectively.
For a $\fh$-module $M$, we denote by $M_\gl$ for $\gl \in \fh^*$ the weight space
\[ M_\gl = \{ v \in M \mid hv = \langle \gl, h\rangle v \ \text{for} \ h \in \fh\}.
\]
We call a vector $v$ of $M$ a \textit{weight vector} if $v \in M_{\gl}$ for some $\gl \in \fh^*$.

Denote by $V(\gl)$ the irreducible highest weight $\fg$-module with highest weight $\gl \in P^+$.
For $w\in W$, let $u_{w\gl}$ be a nonzero vector of the one-dimensional weight space $V(\gl)_{w\gl}$.

\begin{Def} \normalfont
  For $\gl \in P^+$ and $w \in W$, the $\fb$-submodule
  \[ V_w(\gl) = U(\fb) u_{w\gl} \subseteq V(\gl)
  \]
  is called the \textit{Demazure module} associated with $\gl$ and $w$.
\end{Def}

For a subset $J \subseteq I$, we denote by $\fp_J$ the parabolic subalgebra associated with $J$, which is 
the Lie subalgebra of $\fg$ generated by $\fb$ and $\{f_i \mid i \in J\}$.
Note that $\fp_\emptyset = \fb$.
If $J = \{i\}$, we denote $\fp_J$ by $\fp_i$.
In this article, we denote by $\Mod{\fp_J}$ the category of finite-dimensional $\fh$-semisimple left $\fp_J$-modules.
If $\langle w\gl, \ga_i^{\vee}\rangle \le 0$, then the Demazure module $V_w(\gl)$ is preserved by the action of $f_i$.
Hence $V_w(\gl)$ belongs to $\Mod{\fp_J}$ if $\langle w\gl, \ga_i^{\vee} \rangle \le 0$ holds for all $i \in J$.

\begin{Def} \normalfont \label{Def:Demazure_flag}
  Let $M \in \Mod{\fb}$.
  A filtration $0 = M_0 \subseteq M_1 \subseteq \cdots \subseteq M_k = M$ is called a \textit{Demazure flag} of $M$
  if every successive quotient $M_i/ M_{i-1}$ is isomorphic to some Demazure module.
\end{Def}

For every $i \in I$, it is known that the restriction functor $\Mod{\fp_i} \to \Mod{\fb}$ has the left adjoint functor
$\cD_i\colon \Mod{\fb} \to \Mod{\fp_i}$ \cite{MR826100,MR1923198}. 
We often regard $\cD_i$ as a functor from $\Mod{\fb}$ to $\Mod{\fb}$ in the obvious way.

\begin{Prop}[{\cite[Remark 8.1.18]{MR1923198}}]
  Let $w \in W$, and $w = s_{i_k} \cdots s_{i_1}$ be its arbitrary reduced expression. 
  Then the functor $\cD_w = \cD_{i_k} \cdots \cD_{i_1}\colon \Mod{\fb} \to \Mod{\fb}$ does not depend 
  on the choice of the expression.
\end{Prop}

\begin{Def}\normalfont
  The functor $\cD_w\colon \Mod{\fb} \to \Mod{\fb}$ is called the \textit{Joseph functor} associated with $w \in W$.
\end{Def}

Note that for every $i \in I$ and $M \in \Mod{\fb}$, 
the canonical $\fb$-module homomorphism $M \to \cD_i M$ is defined as the image of the identity under the bijection
$\Hom_{\Mod{\fp_i}}(\cD_i M, \cD_i M) \stackrel{\sim}{\to} \Hom_{\Mod{\fb}}(M, \cD_i M)$.
The following lemma was proved in \cite[Subsection 2.7 and Lemma 2.8 (iv)]{MR826100} 
for finite-dimensional $\fg$, and the proof goes without any change for general $\fg$:

\begin{Lem}\label{Lem:Elementary_for_Demazure}
  Let $i \in I$.\\
  {\normalfont (i)} For every $M \in \Mod{\fp_i}$, 
      the canonical $\fb$-module homomorphism $M \to \cD_i M$ is an isomorphism. 
      In particular, we have $\cD_i^2 N \cong \cD_i N$
      for every $N \in \Mod{\fb}$. \\
  {\normalfont (ii)} Assume that $0 \to M_1 \to M_2 \to M_3 \to 0$
      is an exact sequence of objects of $\Mod{\fb}$ and $M_3$ is isomorphic to some $\fb$-submodule of a finite-dimensional 
      $\fp_i$-module.
      Then the sequence $0 \to \cD_iM_1 \to \cD_iM_2 \to \cD_iM_3 \to 0$ is exact.
\end{Lem}

For $\gl \in \fh^*$, we denote by $\C_{\gl}$  
the one-dimensional $\fb$-module spanned by a weight vector with weight $\gl$ on which $e_i$ ($i \in I$) acts trivially. 
Note that we have $V_{\id}(\gl) = \C_{\gl}$ for $\gl \in P^+$.
Now we recall the following theorem:

\begin{Thm}[{\cite[Proposition 8.1.17 and Corollary 8.1.26]{MR1923198}}] \label{Thm:construction}
 For every $\gl \in P^+$ and $w \in W$, we have
 \[ \cD_w \C_\gl \cong V_w(\gl)
 \]
 as $\fb$-modules.
\end{Thm}

\begin{Cor} \label{Cor:Demazure_flag}
  Assume that $M \in \Mod{\fb}$ has a Demazure flag.
  Then $\cD_w M$ has a Demazure flag for every $w \in W$.
\end{Cor}

\begin{proof}
  It is enough to show the assertion for $w = s_i$.
  Let $0= M_0 \subseteq M_1 \subseteq \cdots \subseteq M_k = M$ be a Demazure flag of $M$.
  We show the assertion by the induction on $k$.
  For a Demazure module $V_w(\gl)$, we have from Lemma \ref{Lem:Elementary_for_Demazure} (i) 
  and Theorem \ref{Thm:construction} that
  \begin{equation}\label{eq:Demazure}
    \cD_iV_w(\gl) \cong \begin{cases} V_{s_iw}(\gl) & \text{if} \ \ell(s_i w) = \ell(w) +1, \\
                                       V_w(\gl) & \text{if} \ \ell(s_iw) = \ell(w) -1,
                         \end{cases}
  \end{equation}
  where $\ell$ denotes the length function. Hence the assertion for $k =1$ follows.
  Assume $k >1$ and write $M_k / M_{k-1}\cong V_{w_k} (\gl_k) $.
  Recall that $V_{w_k}(\gl_k)$ is defined as a $\fb$-submodule of $V(\gl_k)$. 
  Hence this module is also a $\fb$-submodule of the $\fp_i$-module $U(\fp_i) V_{w_k}(\gl_k) \subseteq V(\gl_k)$,
  which is finite-dimensional since $V(\gl_k)$ is integrable.
  Hence the sequence $0 \to \cD_i M_{k-1} \to \cD_i M \to \cD_i V_{w_k}(\gl_k) \to 0$ is 
  exact by Lemma \ref{Lem:Elementary_for_Demazure} (ii), and then the assertion follows from the induction hypothesis. 
\end{proof}

Let $\Z[P]$ denote the group ring of $P$ with basis $e^\gl$ ($\gl \in P$),
and define for $i \in I$ a linear operator $D_i$ on $\Z[P]$ by
\begin{equation*}
   D_i(f) = \frac{f - e^{-\ga_i}\cdot s_i(f)}{1-e^{-\ga_i}},
\end{equation*}
where $s_i$ acts on $\Z[P]$ by $s_i(e^{\gl}) = e^{s_i(\gl)}$.
We call $D_i$ the \textit{Demazure operator} associated with $i$.
The following lemma is proved by elementary calculations such as $s_i(fg) = s_i(f)s_i(g)$.

\begin{Lem}\label{Lem:invariance}
  Assume that $f \in \Z[P]$ is $s_i$-invariant. 
  Then we have 
  \[ D_i(f g) = f D_i(g) \ \ \ \text{for every} \ g \in \Z[P].
  \]
  In particular, we have $D_i^2 = D_i$.
\end{Lem}

For a finite-dimensional semisimple $\fh$-module $M$ such that $\{\gl \in \fh^*\mid M_\gl \neq 0 \} \subseteq P$, 
we denote the character of $M$ by
\[ \ch M = \sum_{\gl \in P} \dim M_\gl\cdot e^\gl \in \Z[P].
\]
For every reduced expression $w = s_{i_k} \cdots s_{i_1}$ of $w \in W$, 
the operator $D_w = D_{i_k} \cdots D_{i_1}$ on $\Z[P]$ is independent
of the choice of the expression \cite[Corollary 8.2.10]{MR1923198}, and it is known that the character of
a Demazure module is expressed as follows:

\begin{Thm}[{\cite[Theorem 8.2.9]{MR1923198}}] \label{Thm:Character_formula}
  For every Demazure module $V_w(\gl)$, we have
  \[ \ch V_w(\gl) = D_w (e^{\gl}).
  \]  
\end{Thm}

\begin{Cor}\label{Cor:flag_character}
  Assume that $M \in \Mod{\fb}$ has a Demazure flag.
  Then for every $w \in W$, we have
  \[ \ch \cD_w M = D_w\ch M.
  \]
\end{Cor}

\begin{proof}
  By Corollary \ref{Cor:Demazure_flag}, it is enough to show the assertion for $w = s_i$.
  Let $0= M_0 \subseteq M_1 \subseteq \cdots \subseteq M_k = M$ be a Demazure flag of $M$.
  We show the assertion by the induction on $k$.
  Assume $k =1$, and write $M \cong V_w(\gl)$.
  If $\ell(s_i w) = \ell(w) + 1$, we have
  \begin{equation}\label{eq:Demazure2}
    \ch \cD_i V_w(\gl) = \ch V_{s_i w} (\gl) = D_{s_iw}(e^{\gl}) = D_i \ch V_w(\gl)
  \end{equation}
  by Theorem \ref{Thm:construction} and Theorem \ref{Thm:Character_formula}.
  If $\ell(s_i w) = \ell(w) -1$, we have $\ch \cD_i V_w(\gl) = \ch V_w(\gl)$ by (\ref{eq:Demazure}), and 
  \[ D_i \ch V_w(\gl) = D_i^2 \ch V_{s_iw}(\gl) =D_i \ch V_{s_iw}(\gl) = \ch V_w(\gl)
  \]
  by (\ref{eq:Demazure2}) and Lemma \ref{Lem:invariance}.
  Hence the assertion for $k =1$ follows.
  Assume $k >1$ and write $M_k / M_{k-1} \cong V_{w_k}(\gl_k)$. 
  Then as proved in the proof of Corollary \ref{Cor:Demazure_flag}, 
  the sequence $0 \to \cD_i M_{k-1} \to \cD_i M \to \cD_i V_{w_k}(\gl_k) \to 0$
  is exact. Hence the assertion follows from the induction hypothesis.
\end{proof}

The following theorem is obtained by taking the classical limit of \cite[Theorem 5.22]{MR2214249}:

\begin{Thm}\label{Thm:Joseph}
  Assume that $\fg$ is symmetric {\normalfont (}i.e.\ the Cartan matrix of $\fg$ is symmetric{\normalfont )}. 
  Then for every $\gl, \mu \in P^+$ and $w \in W$,
  the $\fb$-module $\C_{\gl} \otimes V_w(\mu)$ has a Demazure flag.
\end{Thm}

\begin{Rem} \normalfont
  In \cite[Theorem 5.22]{MR2214249}, a given Kac-Moody Lie algebra is assumed to be simply laced. 
  This assumption, however, is used only in \cite[Lemma 3.14]{MR2214249} to apply 
  a positivity result of Lusztig,
  and we can check that the proof of this positivity result in \cite[\S 22.1.7]{MR1227098} goes without any change for 
  every symmetric Kac-Moody Lie algebra.
  Hence \cite[Theorem 5.22]{MR2214249} holds for symmetric Kac-Moody Lie algebras,
  and so does the above theorem. 
\end{Rem}

From the theorem, we have the following corollary since the functor $\C_\gl \otimes -$ is exact:

\begin{Cor}\label{Cor:Joseph}
  Assume that $\fg$ is symmetric.
  If $M \in \Mod{\fb}$ has a Demazure flag, then $\C_{\gl} \otimes M$ for $\gl \in P^+$ has a Demazure flag.
\end{Cor}

\section{Nontwisted affine Lie algebra}\label{Nontwisted affine Lie algebra}

From this section, we assume that $\fg$ is a nontwisted affine Lie algebra with $I= \{0,1,\ldots,n\}$ unless stated otherwise.
Let $\gG$ be the Dynkin diagram of $\fg$, $A = (a_{ij})_{i,j \in I}$ the Cartan matrix of $\fg$, $\gD \subseteq \fh^*$ 
the root system of $\fg$, and $\gD^+ \subseteq \gD$ the set of positive roots corresponding to $\fb$.
In this article, we use the Kac's labeling of nodes of $\gG$ in \cite[Section 4.8]{MR1104219}.
Let $(a_0,\dots,a_n)$ (resp. $(a_0^{\vee},\dots,a_n^{\vee})$) be the unique sequence of relatively prime positive integers 
satisfying 
\[ \sum_{j \in I} a_{ij}a_j = 0 \ \ \ \text{for all $i \in I$} \ (\text{resp.} 
   \sum_{i \in I} a_i^{\vee} a_{ij} = 0 \ \ \ \text{for all $j \in I$}).
\]
Let $d \in \fh$ be the degree operator, which is any element satisfying 
$\langle \ga_i, d \rangle = \gd_{0i}$ for $i \in I$,
$K = \sum_{i \in I} a_i^{\vee}\ga_i^{\vee} \in \fh$ the canonical central element,
and $\gd = \sum_{i \in I} a_i \ga_i \in \fh^*$ the null root.
For each $i \in I$, let $\gL_i \in P^+$ be the fundamental weight, which satisfies 
\[ \langle \gL_i, \ga_j^{\vee} \rangle = \gd_{ij} \ \ \text{for} \ j \in I \ \ \ \text{and} \ \ \ 
   \langle \gL_i, d \rangle = 0.
\]
Note that we have
\[ P = \sum_{i \in I} \Z \gL_i + \C \gd \ \ \ \text{and} \ \ \ P^+ = \sum_{i \in I} \Z_{\ge 0} \gL_i + \C \gd.
\]
Let $( \ , \ )$ be the $W$-invariant symmetric bilinear form on $\fh^*$ defined by
\[ (\ga_i, \ga_j) = a_i^{\vee}a_i^{-1}a_{ij}, \ \ \ (\ga_i, \gL_0) =\gd_{0i} \ \ \text{for}\ i,j\in I \ \ \ \text{and} \ \ \ 
   (\gL_0, \gL_0) = 0.
\]
We denote by $\fn^+$ (resp.\ $\fn^-$) the Lie subalgebra of $\fg$
generated by $\{e_i \mid i \in I\}$ (resp.\ $\{f_i \mid i \in I\}$),
and by $\fg_\ga$ for $\ga \in \gD$ the root space of $\fg$.

Let $I_0 = I \setminus \{ 0 \}$, and $\fg_0 \subseteq \fg$ be the simple Lie subalgebra generated 
by $\{e_i,f_i\mid i \in I_0\}$
with Cartan subalgebra $\fh_0 \subseteq \fh$ and Weyl group $W_0 \subseteq W$.
Let $\gD_0 \subseteq \gD$ be the root system of $\fg_0$, $\gD_0^+ = \gD_0 \cap \gD^+$ the set of positive roots,
$P_0 \subseteq \fh^*_0$ the weight lattice of $\fg_0$, $P_0^+\subseteq P_0$ the set of dominant integral weights,
$Q_0 \subseteq P_0$ the root lattice of $\fg_0$, and $Q_0^+ = \sum_{i \in I_0} \Z_{\ge 0} \ga_i$.
Denote by $\varpi_i \in \fh^*_0$ and $\varpi_i^{\vee}\in\fh_0$ ($i \in I_0$) 
the fundamental weights and the fundamental coweights of $\fg_0$ respectively. 
For the notational convenience, we set $\varpi_0 = 0$ and $\varpi_0^{\vee} = 0$.
We often regard $\fh^*_0$ as a subspace of $\fh^*$ by setting $\langle \fh_0^*, K\rangle =\langle \fh_0^*$, 
$d\rangle = 0$.
Then we have $\fh^* = \fh_0^* \oplus \C \gL_0 \oplus \C \gd$.
We denote by $\gt = \sum_{i \in I_0} a_i \ga_i$ the highest root of $\fg_0$,
and by $w_0$ the longest element of $W_0$.
Let $\fn^+_0$ (resp.\ $\fn^-_0$) be the Lie subalgebra of $\fg_0$ generated by $\{ e_i \mid i \in I_0\}$ 
(resp.\ $\{f_i \mid i \in I_0\}$).
Set $e_{\ga_i} = e_i$ and $e_{-\ga_i} = f_i$ for $i \in I_0$, and 
for each $\ga \in \gD_0^+ \setminus \{\ga_i\mid i \in I_0\}$ fix nonzero vectors $e_{\pm\ga} \in \fg_{\pm \ga}$, 
$\ga^{\vee} \in \fh_0$ so that
\[ [e_\ga, e_{-\ga}] = \ga^{\vee}, \ \ \ [\ga^{\vee}, e_{\pm \ga}] = \pm 2 e_{\pm\ga}.
\]
Recall that there exists a unique Lie algebra isomorphism 
\[ \fg \stackrel{\sim}{\to} \fg_0 \otimes \C [t,t^{-1}] \oplus \C K \oplus \C d
\]
satisfying
\begin{equation*} e_i \mapsto e_i \otimes 1, \ f_i\mapsto f_i \otimes 1 \ (i \in I_0),\
               e_0 \mapsto e_{-\gt}\otimes t, \ f_0 \mapsto e_\gt \otimes t^{-1},\ d \mapsto d. 
\end{equation*}
The Lie algebra structure of $\fg_0 \otimes \C [t,t^{-1}] \oplus \C K \oplus \C d$ is defined by
\begin{align*} [x \otimes t^m + & a_1 K + b_1 d, y \otimes t^n + a_2 K + b_2 d] \\
               &= [x, y] \otimes t^{m+n} + n b_1 y \otimes t^{n} - m b_2 x \otimes t^{m} + m\gd_{m, -n} (x,y)K.
\end{align*}
In the sequel, we always identify these two Lie algebras via the above isomorphism.
It should be noted that we have
\[ \fp_{I_0} = \fg_0 \otimes \C[t] \oplus \C K \oplus \C d.
\]
Set $\fg' = [\fg, \fg]$,
\[ \fp_J' = \fp_J \cap \fg' \ \ \text{for} \ J \subseteq I, \ \ \ \text{and} \ \ \ \fh' = \fh \cap \fg'.
\]
Note that we have $\fp_J = \fp_J' \oplus \C d$ and $\fh' = \fh_0 \oplus \C K$.
Let $\cl \colon \fh^* \to (\fh')^* = \fh^*/ \C \gd$ denote the canonical projection, and set $P_{\cl} = \cl(P)$.
Since $W$ fixes $\gd$, $W$ acts on $\fh^*/\C \gd$ and $P_{\cl}$.
For $\ell \in \Z$, we denote by $P_\cl^\ell$ the subset $\{\gl \in P_{\cl} \mid \langle \gl, 
K \rangle = \ell \}$ of $P_{\cl}$.

As \cite[(6.5.2)]{MR1104219}, we define for $\gl \in P_0$ an endomorphism $t_{\gl}$ of $\fh^*$ by
\begin{equation} \label{eq:translation}
  t_{\gl}(\mu) = \mu + \langle \mu, K \rangle \gl - \big((\mu, \gl) + 
  \frac{1}{2}(\gl, \gl)\langle  \mu, K \rangle\big)\gd.
\end{equation}
The map $\gl \mapsto t_{\gl}$ defines an injective group homomorphism from $P_0$ to the group of 
linear automorphisms of $\fh^*$ orthogonal with respect to $( \ , \ )$.
Let $c_i = a_i/a_i^{\vee}$ for $i \in I_0$,
and define the sublattices $M$ and $\wti{M}$ of $P_0$ by
\[ M = \sum_{w \in W_0} \Z w (\gt), \ \ \ \widetilde{M} = \bigoplus_{i \in I_0} \Z c_i \varpi_i.
\]
Let $T(M)$ and $T(\wti{M})$ be the subgroups of $\mathrm{GL}(\fh^*)$ defined by 
\[ T(M) = \{ t_{\gl} \mid \gl \in M\}, \ \ \ T(\wti{M}) = \{ t_{\gl} \mid \gl \in \wti{M} \} . 
\]
It is known that $W \cong W_0 \ltimes T(M)$ \cite[Proposition 6.5]{MR1104219}.
Define the subgroup $\widetilde{W}$ of $\mathrm{GL}(\fh^*)$ by
\[ \widetilde{W} = W_0 \ltimes T(\widetilde{M}),
\]
which is called the \textit{extended affine Weyl group}.
The action of $\wti{W}$ preserves $\gD$, and $w \in W_0$ and $\gl \in \widetilde{M}$ satisfy
\[ wt_{\gl} w^{-1} = t_{w(\gl)}.
\]
By $\Aut (\gG)$ we denote the group of automorphisms of the Dynkin diagram $\gG$, 
that is, the group of permutations $\tau$ of $I$ satisfying $a_{ij} = a_{\tau(i) \tau(j)}$ for all $i,j \in I$.
Let $\fh_\R^* = \sum_{i \in I} \R \ga_i + \R \gL_0 \subseteq \fh^*$, $C = \{ \gl \in \fh^*_\R \mid ( \gl, \ga_i )  \ge 0 \ 
\text{for all} \ i \in I \}$ be the fundamental chamber,
and $\gS \subseteq \wti{W}$ the subgroup consisting of elements preserving $C$.
Then we have
\[ \wti{W} \cong W \rtimes \gS.
\]
Since $\tau \in \gS$ preserves the set of simple roots,
$\tau$ induces a permutation of $I$ (also denoted by $\tau$) by $\tau(\ga_i) = \ga_{\tau(i)}$,
which belongs to $\Aut(\gG)$ since $( \ , \ )$ is $\tau$-invariant.
By abuse of notation, we denote by $\gS$ both the subgroups of $\wti{W}$ and $\Aut(\gG)$.

\begin{Lem}\label{Lem:transform}
  Let $\tau$ be an arbitrary element of $\gS$, and $\bar{\tau}$ the unique element 
  of $W_0$ such that $\tau \in\bar{\tau}\cdot T(\wti{M})$.
  Then we have
  \[ \tau(\gl + a \gd) = \bar{\tau}(\gl) + \big(a +  \langle \gl, \varpi^\vee_{\tau^{-1}(0)}\rangle\big)\gd \ \ \ 
     \text{for}\ \gl \in P_0 \ \text{and} \ a \in \C.
  \]
\end{Lem}

\begin{proof}
  As $t_{\mu}(\ga_i) \equiv \ga_i$ mod $\Z\gd$
  for every $i\in I$ and $\mu \in \wti{M}$, we have $\tau(\ga_{i}) \equiv \bar{\tau}(\ga_i)$ mod $\Z\gd$.
  This forces $\tau(\ga_i) = \bar{\tau}(\ga_i) + \gd_{\tau(i), 0}\gd$ for $i \in I_0$
  since $\tau$ preserves $\{\ga_0,\ldots,\ga_n\}$ and $\bar{\tau} \in W_0$.
  Now the assertion follows since 
  \begin{align*}
    \tau(\gl + a \gd) &= \sum_{i \in I_0} \langle \gl, \varpi_i^{\vee} \rangle \tau(\ga_{i}) + a \gd \\
                      &= \sum_{i \in I_0} \langle \gl, \varpi_i^{\vee} \rangle \bar{\tau}(\ga_i) + 
                         \big(a + \langle \gl, \varpi^\vee_{\tau^{-1}(0)}\rangle\big) \gd \\
                      & = \bar{\tau}(\gl) + \big(a + \langle \gl, \varpi^\vee_{\tau^{-1}(0)}\rangle\big)\gd.
  \end{align*}
\end{proof}

We define an action of $\gS$ on $\fg$ by letting $\tau \in \gS$ act as a Lie algebra automorphism defined by
\[ \tau(e_i) = e_{\tau(i)}, \ \ \ \tau(\ga_i^{\vee}) = \ga_{\tau(i)}^{\vee}, \ \ \ \tau(f_i) = f_{\tau(i)} \ \ \
   \text{and} \ \ \ \tau(d) = d + \varpi_{\tau(0)}^{\vee}.
\]
This action obviously preserves $\fb$ and $\fb'$.
For a module $M$, we denote by $\tau^*M$ the pull-back of $M$ with respect to $\tau$. 
The image of $v \in M$ under the canonical linear isomorphism $M \to \tau^*M$ is denoted by $\tau^* v$.
Note that this isomorphism maps the weight space $M_{\gl}$ onto $(\tau^*M)_{\tau^{-1}(\gl)}$.

We prepare some notation here.
Let $w \in \wti{W}$ be an arbitrary element, and take unique elements $w' \in W$ and $\tau\in \gS$ so that 
$w = w'\tau$.
Then define a functor $\cD_{w}\colon \Mod{\fb} \to \Mod{\fb}$ by
\[ \cD_{w} = \cD_{w'} \circ \big(\tau^{-1}\big)^*,
\]
and a linear operator $D_w\colon \Z[P] \to \Z[P]$ by
\[ D_w = D_{w'} \circ \tau,
\]
where $\tau$ acts on $\Z[P]$ by $\tau(e^{\gl}) = e^{\tau(\gl)}$.
We set 
\[ V_{w}(\gL) = V_{w'}\big(\tau(\gL)\big) \ \ \ \text{for} \ \gL \in P^+.
\]
Then we have the following lemma:

\begin{Lem}\label{Lem:affine_version}
  {\normalfont(i)}  
    Let $\tau \in \gS$, $\gL \in P^+$ and $w \in W$. Then we have
    \[ \big(\tau^{-1}\big)^* V_w(\gL) \cong V_{\tau w}(\gL).
    \]
  {\normalfont(ii)} 
    Assume $M \in \Mod{\fb}$ has a Demazure flag.
    Then for every $w \in \wti{W}$, $\cD_w M$ has a Demazure flag and we have
    \[ \ch \cD_w M = D_w \ch M.
    \]
\end{Lem}

\begin{proof}
  Since $\big(\tau^{-1}\big)^*V(\gL)$ is an integrable highest weight $\fg$-module with highest weight $\tau(\gL)$,
  we have $\big(\tau^{-1}\big)^*V(\gL) \cong V\big(\tau(\gL)\big)$.
  Moreover, we see that $\big(\tau^{-1}\big)^*V_w(\gL)$ is the $\fb$-submodule of $\big(\tau^{-1}\big)^* V(\gL)$ 
  generated by the weight space with weight $\tau w(\gL)$.
  Hence the assertion (i) follows.
  Then, since the functor $\big(\tau^{-1}\big)^*$ is exact, 
  the first assertion of (ii) follows from (i) and Corollary \ref{Cor:Demazure_flag}.
  Since we have $\ch \big(\tau^{-1}\big)^* M = \tau\, \ch M$, 
  the second one follows from Corollary \ref{Cor:flag_character}. 
\end{proof}

Let $\Z[P_{\cl}]$ denote the group algebra of $P_{\cl}$ with basis $e^{\gl}$ ($\gl \in P_{\cl}$),
and by $\cl$ we also denote the projection from $\Z[P]$ to $\Z[P_{\cl}]$ defined by $e^{\gl} \mapsto e^{\cl(\gl)}$.
For each $w \in W$, a linear operator $\ol{D}_w$ on $\Z[P_{\cl}]$ is defined by 
\[ \cl \circ D_w = \ol{D}_w \circ \cl.
\]

\begin{Lem}\label{Lem:final_lem}
  Let $f \in \Z[P_{\cl}^0]$, and assume $f$ is $W_0$-invariant.
  Then we have 
  \[ \ol{D}_w(f g) = f \ol{D}_w(g) \ \ \ \text{for every} \ g \in \Z[P_{\cl}] \ \text{and} \ w \in W.
  \]
\end{Lem}

\begin{proof}
  Since $f \in \Z[P_{\cl}^0]$, we have $s_0(f) = s_{\gt}(f) =f$, where $s_{\gt}$ denotes the reflection 
  associated with $\gt$.
  Hence $f$ is $s_i$-invariant for all $i \in I$, and then the assertion is easily proved from 
  Lemma \ref{Lem:invariance}.
\end{proof}

We need the following elementary lemmas later:

\begin{Lem}\label{Lem:trivial_lem}
  Let $M \in \Mod{\fb}$.
  For every $w \in \wti{W}$ and $C \in \C$ , we have 
  \[ \cD_w(\C_{C\gd} \otimes M) \cong \C_{C\gd} \otimes \cD_wM.
  \]
\end{Lem}

\begin{proof}
  It is enough to show the assertions for $w = \tau \in \gS$ and $w =s_i$.
  The first case is obvious, and the second one follows since we have
  \begin{align*} 
    &\Hom_{\Mod{\fb}}(\C_{C \gd} \otimes M, M')\cong \Hom_{\Mod{\fb}}(M, \C_{-C \gd} \otimes M')  \\
    &\cong \Hom_{\Mod{\fp_i}}(\cD_iM, \C_{-C\gd} \otimes M')\cong \Hom_{\Mod{\fp_i}}(\C_{C\gd} \otimes \cD_iM, M') 
  \end{align*} 
  for every $M' \in \Mod{\fp_i}$.
\end{proof}

\begin{Lem}\label{Lem:trivial_lem2}
  Let $M_1$ be an object of $\Mod{\fb}$ which is generated by a weight vector $v_1$, $M_2$ an object of $\Mod{\fb}$,
  and $\Phi$ a homomorphism of $\fb'$-modules from $M_1$ to $M_2$ which maps $v_1$ to a weight vector $v_2$.\\
  {\normalfont (i)} For some $C \in \C$, $\Phi$ induces a homomorphism of $\fb$-modules
      from $\C_{C\gd} \otimes M_1$ to $M_2$. \\
  {\normalfont(ii)} Assume further that $M_2$ extends to an object of $\Mod{\fp_i}$ for some $i \in I$.
       Then there exists a homomorphism $\wti{\Phi}$ of $\fp_i'$-modules from $\cD_i M_1$ to $M_2$ satisfying
       $\Phi = \wti{\Phi} \circ \iota$, where $\iota\colon M_1 \to \cD_i M_1$ is the canonical homomorphism of $\fb$-modules. 
\end{Lem}

\begin{proof}
  Let $\gl_1,\gl_2 \in \fh^*$ be the respective weights of $v_1, v_2$,
  and $C = \langle \gl_2 - \gl_1, d \rangle$.
  It is easily seen that the induced map $\Phi'\colon\C_{C\gd} \otimes M_1 \to M_2$ preserves the weights, 
  and hence it is a homomorphism of $\fb$-modules.
  The assertion (i) is proved.
  Then under the assumption of (ii), there exists a homomorphism 
  $\wti{\Phi}'\colon \cD_i(\C_{C\gd} \otimes M_1) \to M_2$ of $\fp_i$-modules
  such that $\Phi' = \wti{\Phi}' \circ \iota$ since $\cD_i$ is left adjoint to the restriction functor 
  $\Mod{\fp_i} \to \Mod{\fb}$.  
  Since $\cD_i(\C_{C\gd} \otimes M_1) \cong \C_{C\gd}\otimes \cD_iM_1$ holds by Lemma \ref{Lem:trivial_lem}, 
  required $\wti{\Phi}$ is obtained by restricting $\wti{\Phi}'$ to $\fp_i'$. 
\end{proof}

\section{Kirillov-Reshetikhin modules}\label{Section:KR}

Following \cite{MR2238884}, we define the following $\fp_{I_0}$-modules:

\begin{Def}\normalfont\label{Def:KR_module}
  For $r \in I_0$ and $\ell \in \Z_{> 0}$, let $KR^{r,\ell}$ be the 
  $\fp_{I_0}$-module generated by a nonzero vector 
  $v_{r,\ell}$ with relations
  \begin{align*}
    (\fn_0^+ \otimes \C[t]) v_{r,\ell} &= 0, \ \ \ (\fh_0 \otimes t\C[t])v_{r,\ell} = 0, \ \ \
    hv_{r,\ell} = \langle \ell \varpi_r, 
    h \rangle v_{r,\ell} \ \ \text{for} \ h \in \fh, \\
    f_r^{\ell+1}v_{r,\ell} &= (f_r\otimes t)v_{r,\ell} = 0 \ \ \
    \text{and} \ \ \ f_{i}v_{r,\ell} = 0 \ \ \text{for} \ i \in I_0 \setminus \{r\}.    
  \end{align*}
  We call $KR^{r,\ell}$ the \textit{Kirillov-Reshetikhin module} (KR module for short) for $\fp_{I_0}$ 
  associated with $r$ and $\ell$.
\end{Def}

\begin{Rem}\normalfont
  (i) Kirillov-Reshetikhin modules were originally defined in \cite{MR2238884} as ($\fg_0\otimes \C[t]$)-modules.
      Since we would like to consider $KR^{r,\ell}$ as a $\fb$-module later, 
      we adopt the above definition in this article.
      It is obvious that the restriction of $KR^{r,\ell}$ to $\fg_0\otimes \C[t]$ 
      coincides with the original one.\\
  (ii) 
    As the name indicates, $KR^{r,\ell}$ has strong connections with the Kirillov-Reshetikhin module $W^{r,\ell}$
    for the quantum affine algebra $U_q'(\fg)$ \cite{MR1836791,MR2238884,MR2428305}. 
    We return to this topic in Section \ref{X=M_conjecture}.
\end{Rem}

It is easily seen that the $\fg_0$-submodule $U(\fg_0)v_{r,\ell} \subseteq KR^{r,\ell}$ 
is isomorphic to the irreducible module with highest weight $\ell\varpi_r$.
For each $w \in W_0 \setminus \{ \id \}$, take and fix a nonzero vector  $v_{r,\ell}^w$ 
of this $\fg_0$-submodule whose weight is $w(\ell\varpi_r)$,
and set $v_{r,\ell}^{\id} = v_{r,\ell}$.

\begin{Lem}\label{Lem:KR}
  {\normalfont(i)} Let $w \in W_0$. For $\ga \in Q^+_0$ and $k \in \Z$, we have
  \[ KR^{r,\ell}_{\ell w(\varpi_r) + w(\ga) + k \gd} =
                                           \begin{cases}  \C v_{r,\ell}^{w} & \text{if} \ \ga = 0 \ \text{and} \ k = 0, \\
                                                          \{ 0 \}           & \text{otherwise}.
                                           \end{cases}
  \]
  {\normalfont(ii)} As a $\fb$-module, $KR^{r,\ell}$ is generated by $v_{r,\ell}^{w_0}$.
\end{Lem}

\begin{proof}
  The assertion (i) for $w = \id$ follows by definition.
  Since $KR^{r,\ell}$ is finite-dimensional \cite[Theorem 1.2.2]{MR2271991}, 
  its weight set and dimensions of weight spaces are $W_0$-invariant.
  From this, (i) for general $w$ is proved.
  Then (i) for $w = w_0$ implies (ii).
\end{proof}

The following theorem, which easily follows from \cite[Theorem 4]{MR2323538}, connects KR modules with Demazure modules.

\begin{Thm}\label{Thm:KR-Demazure}
  Let $r$ be an arbitrary element of $I_0$ and $\ell$ a positive integer satisfying $c_r^{-1} \ell \in \Z$.
  Then there exists an isomorphism
  \[ \C_{c_r^{-1}\ell\gL_0} \otimes KR^{r,\ell} \stackrel{\sim}{\to} V_{t_{c_rw_0(\varpi_r)}}(c_r^{-1}\ell\gL_0)
  \]
  of $\fp_{I_0}'$-modules which maps $u_{c_r^{-1}\ell\gL_0} \otimes v^{w_0}_{r,\ell}$ to 
  $u_{t_{c_rw_0(\varpi_r)}(c_r^{-1}\ell\gL_0)}$.
\end{Thm}

\begin{Rem}\normalfont
  If $\fg_0$ is of $ADE$ type, then $c_r = 1$ holds for all $r \in I_0$.
  Hence the above isomorphism follows for all KR modules in this case.
\end{Rem}

\section{Fusion product}\label{Fusion product}

The fusion product was defined in \cite{MR1729359} as a construction of a graded cyclic ($\fg_0 \otimes \C[t]$)-module.
Here we slightly reformulate it as a construction of a $\fp_{I_0}$-module.
Note that $U(\fp_{I_0}')$ has a natural grading defined by
\begin{equation*}
  U(\fp_{I_0}')^k = \{ X \in U(\fp_{I_0}') \mid [d, X] = kX \},
\end{equation*}
from which we define a natural filtration on $U(\fp_{I_0}')$ by
\[ U(\fp_{I_0}')^{\le k} = \bigoplus_{q \le k} U(\fp_{I_0}')^q.
\]
Let $M$ be a cyclic $\fp_{I_0}'$-module with a generator $v$,
and denote by $F^k_v(M)$ for $k \in \Z_{\ge -1}$ the subspace $U(\fp_{I_0}')^{\le k}v$ of $M$.
(Note that $F^{-1}_v(M) = 0$.)
Then the associated $\fp_{I_0}$-module $\mathrm{gr}_{F_v}(M)$ is defined by
\[ \mathrm{gr}_{F_v}(M) = \bigoplus_{k \ge 0} F^k_v(M) / F^{k-1}_v(M),
\]
where $d$ acts on $F_v^k(M)/F_v^{k-1}(M)$ as multiplication by $k$.

Now we recall the definition of fusion products. 
Let $M^1, \dots, M^p$ be a sequence of objects of $\Mod{\fp_{I_0}}$ such that each $M^j$ is generated by a weight vector $v_j$,
and $c_1, \dots, c_p$ pairwise distinct complex numbers.
For each $1 \le j \le p$, define $\fp_{I_0}'$-module $M_{c_j}^j$ by the pullback $\varphi_{c_j}^*M^j$, where 
$\varphi_{c}$ is an automorphism of $\fp_{I_0}'$ defined by
\[ x \otimes t^k \mapsto x \otimes (t +c)^k \ \text{for} \ x \in \fg_0, \ \ \   K \mapsto K.
\]
For $w \in M^j$, we denote by $w'$ its image under the canonical map $M^j \to M^j_{c_j}$.
As shown in \cite[Proposition 1.4]{MR1729359}, $M^1_{c_1} \otimes \dots \otimes M^p_{c_p}$ is 
a cyclic $\fp_{I_0}'$-module generated by $v_1' \otimes \dots \otimes v_p'$,
and we define a $\fp_{I_0}$-module $M^1_{c_1,v_1} * \dots * M^p_{c_p,v_p}$ by
  \[ M^1_{c_1,v_1} * \dots * M^p_{c_p,v_p} = \mathrm{gr}_{F_{v_1'\otimes \cdots \otimes v_p'}}
     (M^1_{c_1} \otimes \dots \otimes M^p_{c_p}).
  \]
When the parameters or generators are clear from the context,
we write simply as $M_{c_1}^1 * \cdots *M_{c_p}^p$, $M_{v_1}^1*\cdots *M_{v_p}^p$ or $M^1 * \cdots * M^p$.

\begin{Def}[\cite{MR1729359}]\normalfont 
  The $\fp_{I_0}$-module $M^1 * \dots * M^p$ is called the \textit{fusion product}.
\end{Def}

\begin{Lem}\label{Lem:elementary}
  {\normalfont(i)} As $(\fg_0\oplus \C K)$-modules, 
    \[ M^1* \cdots *M^p \cong M^1 \otimes \cdots \otimes M^p.
    \]
  {\normalfont(ii)} Let $w_1,\ldots, w_p$ be weight vectors of $M_1,\ldots,M_p$ respectively, and assume that 
    \[ U(\fg_0)(w_1 \otimes \cdots \otimes w_p) = U(\fg_0)(v_1\otimes\cdots\otimes v_p).
    \]
    Then we have
    \[ M^1_{w_1} * \cdots *M^p_{w_p} \cong M^1_{v_1} * \cdots * M^p_{v_p}.
    \]
  {\normalfont(iii)} For any $c \in \C$, we have 
    \[ M^1_{c_1+c} * \cdots * M^p_{c_p+c} \cong M^1_{c_1} * \cdots * M^p_{c_p}.
    \]  
  {\normalfont(iv)} $M_{c_1,v_1}^1$ {\normalfont(}the fusion product of a single module $M^1${\normalfont)} 
    is isomorphic to $M^1$ as a $\fp_{I_0}'$-module.
\end{Lem}

\begin{proof}
  The assertions (i) and (ii) easily follow from the definition.
  As the automorphism $\varphi_c$ preserves the subspace $U(\fp_{I_0}')^{\le k}$, we have
  \[ \varphi_c\big(U(\fp_{I_0}')^{\le k}\big)(v_1' \otimes \cdots \otimes v_p')
     = U(\fp_{I_0}')^{\le k}(v_1' \otimes \cdots \otimes v_p').
  \]
  Since the left hand side is equal to the filtration of $M^1_{c_1 + c} \otimes \cdots \otimes M^p_{c_p+c}$,
  (iii) is proved.
  When $c_1 = 0$, (iv) obviously follows.
  Then the assertion follows in general by (iii).
\end{proof}

\section{Statement of the main theorem}\label{Section:Main}

Now we state the main theorem of this article.
This is a generalization of \cite[Theorem 2.5]{MR1729359}, in which the case $\fg_0 = \mathfrak{sl}_2$ is proved.

\begin{Thm}\label{Thm:Main_theorem}
  Assume that $\fg_0$ is of $ADE$ type.
  Let $KR^{r_1,\ell_1}, \ldots, KR^{r_p,\ell_p}$ be a sequence of KR modules such that $\ell_1 \le \cdots \le \ell_p$,
  and set $v_j = v_{r_j,\ell_j} \in KR^{r_j,\ell_j}$ for $1 \le j \le p$.
  Then for arbitrary pairwise distinct complex numbers $c_1, \ldots,c_p$, there exists an isomorphism of $\fb$-modules
  \begin{align}\label{eq:Main_theorem}
     &\C_{\ell_p\gL_0 + C \gd}\otimes \left(KR^{r_p,\ell_p}_{c_p,v_p} * \cdots *KR^{r_2,\ell_2}_{c_2,v_2}
      *KR^{r_1,\ell_1}_{c_1,v_1}\right) \\ 
     & \
     \cong\cD_{t_{w_0(\varpi_{r_p})}}\Big(\C_{(\ell_p-\ell_{p-1})\gL_0} \otimes \cdots \otimes \cD_{t_{w_0(\varpi_{r_2})}}
     \big( \C_{(\ell_2-\ell_1)\gL_0} \otimes \cD_{t_{w_0(\varpi_{r_1})}}\C_{\ell_1\gL_0}\big) \cdots \Big) \nonumber
  \end{align}
  with some constant $C$.
\end{Thm}

We postpone the proof of this theorem to the latter part of this article.
We see from Corollary \ref{Cor:Joseph} and  Lemma \ref{Lem:affine_version} (ii)
that the right hand side of (\ref{Thm:Main_theorem}) has a Demazure flag.
Hence we can prove inductively using Lemma \ref{Lem:affine_version} (ii) that the following equation holds:

\begin{Cor}\label{Cor:Main_cor}
  Under the notation and the assumptions of Theorem \ref{Thm:Main_theorem}, we have 
  \begin{align*}
    &e^{\ell_p \gL_0 + C\gd}\ch KR^{r_p,\ell_p} * \cdots * KR^{r_2,\ell_2} * KR^{r_1,\ell_1} \\ 
    &\ \ \ \ \ \ \ = D_{t_{w_0(\varpi_{r_p})}}\Big(e^{(\ell_p-\ell_{p-1})\gL_0} \cdots
    D_{t_{w_0(\varpi_{r_2})}}
     \big( e^{(\ell_2-\ell_1)\gL_0} \cdot D_{t_{w_0(\varpi_{r_1})}}(e^{\ell_1\gL_0})\big) \cdots \Big).
  \end{align*}
\end{Cor}

\begin{Rem}\normalfont
  The right hand side of (\ref{eq:Main_theorem}) also appeared in \cite{MR1887117}.
  In the article, it was proved that this module, which was called a \textit{generalized Demazure module},
  is isomorphic to the space of global sections of a certain line bundle on a \textit{Bott-Samelson variety}.
\end{Rem}

\section{$X=M$ conjecture}\label{X=M_conjecture}

In this section, we give an important application of Theorem \ref{Thm:Main_theorem},
the proof of the $X=M$ conjecture for type $A_n^{(1)}$ and $D_n^{(1)}$.
Here, we assume that $\fg$ is a general (possibly twisted) affine Lie algebra.

For a sequence ${\bm\nu} = \big((r_1,\ell_1),\ldots,(r_p,\ell_p)\big)$ of elements of $I_0 \times \Z_{>0}$ 
and $\mu \in P_0^+$,
denote by $M(\bm{\nu},\mu,q)\in \Z[q^{-1}]$ the associated \textit{fermionic form}
(see \cite{MR1903978,MR1745263} for definition, in which the fermionic form is denoted by $M(W,\mu,q)$
with $W= \bigotimes_{1 \le j \le p}W^{r_j,\ell_j}$).
The most important result concerning fermionic forms in this article is the following theorem proved by 
Di Francesco and Kedem using the result of \cite{MR2290922}:

\begin{Thm}[\cite{MR2428305}]\label{Thm:DFK}
  Assume that $\fg$ is of nontwisted type. 
  For a sequence $\bm{\nu}= \big((r_1,\ell_1),\ldots, (r_p,\ell_p)\big)$ of elements of $I_0 \times \Z_{>0}$
  and pairwise distinct complex numbers $c_1, \ldots, c_p$, we have
  \[ \ch KR^{r_p,\ell_p}_{c_p} * \cdots * KR^{r_2,\ell_2}_{c_2} * KR^{r_1,\ell_1}_{c_1} = 
     \sum_{\mu \in P_0^+} M(\bm{\nu},\mu,q)
     \ch V_{\fg_0}(\mu),
  \]
  where we set $q = e^{-\gd}$ and denote by $V_{\fg_0}(\mu)$ the irreducible $\fg_0$-module with highest weight $\mu$.  
\end{Thm}

\begin{Rem}\normalfont
  In \cite{MR2428305}, the above theorem was proved under the assumption that the dimension of each $KR^{r_j,\ell_j}$
  is equal to that of the corresponding KR-module $W^{r_j,\ell_j}$ for the quantum affine algebra $U_q'(\fg)$.
  We can see from the pentagon of identities in \cite[Subsection 1.2]{Ke} that 
  this assumption in fact holds in general.     
\end{Rem}

From Corollary \ref{Cor:Main_cor} and Theorem \ref{Thm:DFK}, we have the following corollary:

\begin{Cor}\label{Cor:M-side}
  Assume that $\fg$ is of nontwisted type and $\fg_0$ is of $ADE$ type. 
  Let $\bm{\nu}= \big((r_1,\ell_1),\ldots, (r_p,\ell_p)\big)$ be a sequence of elements of 
  $I_0 \times \Z_{>0}$ such that $\ell_1 \le \cdots \le \ell_p$. 
  Then we have
  \begin{align*}
    q^C&e^{\ell_p\gL_0}\sum_{\mu \in P_0^+} M(\bm{\nu}, \mu,q) \ch V_{\fg_0}(\mu)\\ 
    &= D_{t_{w_0(\varpi_{r_p})}}\Big(e^{(\ell_p-\ell_{p-1})\gL_0} \cdots
    D_{t_{w_0(\varpi_{r_2})}}
    \big( e^{(\ell_2-\ell_1)\gL_0} \cdot D_{t_{w_0(\varpi_{r_1})}}(e^{\ell_1\gL_0})\big) \cdots \Big)
  \end{align*}
  with some constant $C$, where we set $q = e^{-\gd}$.
\end{Cor}

Next, we recall the definition of one-dimensional sums.
Denote by $B^{r,\ell}$ the \textit{Kirillov-Reshetikhin crystal} (KR crystal for short) 
associated with $r\in I_0$ and $\ell\in \Z_{>0}$.
For nonexceptional $\fg$, it is known that $B^{r,\ell}$ is perfect if and only if $\ell \in c_r \Z$ \cite{MR2642564}.
Let $B = B^{r_p,\ell_p} \otimes \cdots \otimes B^{r_1,\ell_1}$ be a tensor product of KR crystals, 
and $D = D_B\colon B \to \Z$ the \textit{energy function} defined on $B$. 
For the definitions of these objects, see \cite{MR1903978,MR1745263}.
The \textit{one-dimensional sum} $X(B, \mu, q) \in \Z[q, q^{-1}]$
for $\mu \in P_0^+$ is defined by
\[ X(B, \mu ,q) = \sum_{\begin{smallmatrix} b \in B \\ \wti{e}_i b= 0 \ (i \in I_0) \\ \wt(b) = 
   \mu \end{smallmatrix}} q^{D(b)},
\]
where $\wti{e}_i$ are Kashiwara operators.
In \cite[Corollary 7.3]{Na3}, the following proposition was proved:

\begin{Prop}\label{Prop:X-side}
  Assume that $\fg$ is of nonexceptional type.
  Let $B = B^{r_p,c_{r_p}\ell_p} \otimes \cdots \otimes B^{r_1, c_{r_1}\ell_1}$ be 
  a tensor product of perfect KR crystals 
  such that $\ell_1 \le \cdots \le \ell_p$. 
  Then we have
  \begin{align*}
    q^{C'}e^{\ell_p\gL_0}\sum_{\mu \in P_0^+} X(B, \mu,q) \, &\ch V_{\fg_0}(\mu)
    = D_{t_{c_{r_p}w_0(\varpi_{r_p})}}\Big( e^{(\ell_p-\ell_{p-1})\gL_0} \cdots \\ 
    & \cdots  D_{t_{c_{r_2}w_0(\varpi_{r_2})}}
    \big( e^{(\ell_2 -\ell_1)\gL_0}\cdot D_{t_{c_{r_1}w_0(\varpi_{r_1})}}(e^{\ell_1\gL_0})\big)\cdots \Big)
  \end{align*}
  with some constant $C'$, where we set $q = e^{-\gd}$.
\end{Prop}

Now we show the following theorem from the above results, which is the $X=M$ conjecture 
presented in \cite{MR1903978,MR1745263}.
This result for $D_n^{(1)}$ is new. This has already been proved for $A_n^{(1)}$ in \cite{MR1890195}, 
but our approach is quite different from theirs:

\begin{Thm}\label{Thm:X=M_conj}
  Assume that $\fg$ is of nontwisted, nonexceptional type and $\fg_0$ is of $ADE$ type
  {\normalfont(}i.e. $\fg = A_n^{(1)}$ or $D_n^{(1)}${\normalfont)}.
  Let ${\bm\nu} = \big((r_1,\ell_1),\ldots,(r_p,\ell_p)\big)$ be a sequence of elements of $I_0 \times \Z_{>0}$, 
  and $B= B^{r_p,\ell_p} \otimes \cdots \otimes B^{r_1,\ell_1}$.
  Then for every $\mu \in P_0^+$, we have
  \[ q^{-D(u(B))}X(B,\mu,q) = M(\bm{\nu},\mu,q),
  \]
  where $u(B)$ denotes the unique element of $B$ whose weight is $\sum_{1 \le j \le p} \ell_j\varpi_j$.
\end{Thm}

\begin{proof}
  Recall that energy functions and one-dimensional sums are invariant under reordering of the given sequence by 
  \cite[Proposition 2.15]{MR1973369}, and so are fermionic forms by definition.
  Hence we may assume $\ell_1 \le \cdots \le \ell_p$.
  Then as $c_{r} = 1$ holds for all $r\in I_0$, 
  we have from Corollary \ref{Cor:M-side} and Proposition \ref{Prop:X-side} that
  \begin{equation}\label{eq:identity}
    q^{C''}\sum_{\mu \in P_0^+} X(B, \mu,q) \ch V_{\fg_0}(\mu) = \sum_{\mu \in P_0^+} 
    M(\bm{\nu}, \mu,q) \ch V_{\fg_0}(\mu)
  \end{equation}
  with some constant $C''$, which implies 
  \[ q^{C''}X(B,\mu,q) = M(\bm{\nu},\mu,q)
  \]
  for every $\mu$ since the characters of irreducible $\fg_0$-modules are linearly independent.
  It remains to show $C'' = -D\big(u(B)\big)$.
  Let $\gl = \sum_{1 \le j \le p} \ell_j\varpi_j$.
  Since
  \[ \sum_{\mu \in P_0^+} X(B, \mu,q) \ch V_{\fg_0}(\mu) = \sum_{b \in B} q^{D(b)}e^{\wt(b)}
  \]
  holds by definition, the coefficient of $e^{\gl}$ in the left hand side of (\ref{eq:identity}) is 
  equal to $q^{C'' + D(u(B))}$.
  On the other hand, we easily see from Theorem \ref{Thm:DFK} that the coefficient of 
  $e^{\gl}$ in the right hand side is $1$. Hence we have $q^{C'' + D(u(B))} = 1$, which implies $C'' = -D\big(u(B)\big)$.
  The theorem is proved.
\end{proof}

\section{$\fb$-fusion product}\label{Section:b-fusion}

We devote the rest of this article to prove Theorem \ref{Thm:Main_theorem}.
In this section, we introduce a construction of a $\fb$-module, which we call the $\fb$-fusion product,
defined by modifying the definition of the fusion product in Section \ref{Fusion product}.
This construction is essentially used in the proof of the theorem.
Here we only assume that $\fg$ is a nontwisted affine Lie algebra (that is, $\fg_0$ is allowed to be of type $BCFG$)
since the definition of $\fb$-fusion products makes sense in this setting.

Let $M^1,\ldots,M^p$ be a sequence of objects of $\Mod{\fp_{I_0}}$ such that each $M^j$ is generated (as a $\fp_{I_0}$-module) 
by a weight vector $v_j$, and $N$ an object of $\Mod{\fb}$ which is generated (as a $\fb$-module) by a weight vector $u$.
For pairwise distinct \textit{nonzero} complex numbers $c_1,\ldots,c_p$, 
define a $\fp_{I_0}'$-module $M^1_{c_1} \otimes \cdots \otimes M_{c_p}^p$ as in Section \ref{Fusion product}.

\begin{Lem}\label{Lem:cyclicity}
  $N \otimes M_{c_1}^1 \otimes \cdots \otimes M_{c_p}^p$ is generated by the vector
  $u \otimes v_1' \otimes \cdots \otimes v_p'$ as a $\fb'$-module.
\end{Lem}

\begin{proof}
  The proof is similar to that of \cite[Proposition 1.4]{MR1729359}. 
  Here we give it for completeness.
  Let $z \in N$ and $w_j \in M^j$ for $1 \le j\le p$ be arbitrary vectors, and set $\wti{w} = z \otimes 
  w_1'\otimes \cdots \otimes w_p' \in N \otimes M_{c_1}^1 \otimes \cdots \otimes M_{c_p}^p$
  (recall that $w_j'$ denotes the image of $w_j$ under the canonical map).
  First we show that for every $x \in \fg_0$, $1\le j \le p$ and $q \in \Z_{\ge 0}$, 
  there exists a vector $X_j[q] \in \fb'$ such that
  \[ X_j[q]\wti{w}= z \otimes w_1'\otimes \cdots \otimes \big((x\otimes t^q) w_j\big)' \otimes \cdots \otimes w_p'.
  \]
  Since $N$ and $M^j$ are finite-dimensional and $\fh$-semisimple, there exists a sufficiently large positive integer $L$ 
  such that $\fg_0 \otimes t^L\C[t]$ acts trivially on them.
  Since $X_j[q] = 0$ satisfies the above equation for $q \ge L$, we may assume $q \le L-1$.
  Note that we have for $k \ge 0$ that
  \[ (x \otimes t^{k+L})\wti{w} = \sum_{\begin{smallmatrix}1 \le j \le p \\ 0\le q \le L-1\end{smallmatrix}} 
     \begin{pmatrix} \, k+L \,\\ \,q\, \end{pmatrix} c_j^{k+L-q} \cdot z \otimes w_1'\otimes 
     \cdots \otimes \big((x \otimes t^q)w_j\big)' \otimes \cdots \otimes w_p'.
  \] 
  Consider the matrix of the coefficients:
  \[ \Bigg(\begin{pmatrix} \, k+L \, \\ \,q\, \end{pmatrix} c_j^{k+L-q}  \Bigg)_{k, (j,q)}, \ \ 1 \le k \le pL, 
     \ 1 \le j \le p,\ 0 \le q \le L-1.
  \]
  Since the determinant of this matrix is equal to 
  \[ \prod_{\begin{smallmatrix}1 \le j \le p \\ 0 \le q \le L-1\end{smallmatrix}} \left. \frac{1}{q!}
     \left(\frac{\partial}{\partial x_{pq+j}}\right)^q \det(x_\ell^{k+L})_{1\le k,\ell \le pL} 
     \right|_{x_j = x_{p+j} = \cdots = x_{p(L-1) +j} = c_j},
  \]
  we see that this matrix is invertible from the assumption that $c_1, \ldots, c_p$ are nonzero and pairwise distinct.
  Hence $X_j[q]$ can be obtained as a linear combination of $x \otimes t^{k+L}$ ($1 \le k \le pL$), 
  and our assertion is proved.
  From this, we can also prove that for every $y \in \fb'$, 
  there exists a vector $Y \in \fb'$ satisfying
  \[ Y\wti{w}= yz \otimes w_1'\otimes \cdots  \otimes w_p'.
  \]
  Now the lemma obviously follows. 
\end{proof}

Define the subspace $U(\fb')^{\le k} \subseteq U(\fb')$ for $k \in \Z$ similarly as $U(\fp_{I_0}')^{\le k}$.
By considering the filtration 
\begin{equation}\label{eq:filtration}
   \wti{F}_{u \otimes v_1' \otimes \cdots \otimes v_p'}^k(N \otimes M_{c_1}^1 \otimes \cdots \otimes M_{c_p}^p) =
   U(\fb')^{\le k}(u \otimes v_1' \otimes \cdots \otimes v_p')
\end{equation}
of $N \otimes M_{c_1}^1 \otimes \cdots \otimes M_{c_p}^p$, 
we define $\left[N_u * M_{c_1,v_1}^1 * \cdots * M_{c_p,v_p}^p\right]_{\fb}$ as the associated $\fb$-module
\[ \left[N_u * M_{c_1,v_1}^1 * \cdots * M_{c_p,v_p}^p\right]_{\fb} 
   = \mathrm{gr}_{\wti{F}_{u \otimes v_1'\otimes \cdots \otimes v_p'}}(N\otimes M_{c_1}^1 \otimes \cdots \otimes M_{c_p}^p),
\]
which we call the \textit{$\fb$-fusion product}.
We sometimes omit the parameters or the generators when they are clear from the context.
It should be noted that in the definition of the $\fb$-fusion product, only the leftmost module is allowed 
to be a $\fb$-module, and the others are assumed to be $\fp_{I_0}$-modules.
By definition, it is easily seen for every $\ell \in \C$ that
\begin{equation}\label{eq:level}
   \C_{\ell \gL_0} \otimes \left[N * M^1 * \cdots * M^p\right]_{\fb} \cong 
   \left[(\C_{\ell \gL_0} \otimes N) * M^1 * \cdots * M^p\right]_{\fb}.
\end{equation}

In some special cases, an original fusion product is connected to a certain $\fb$-fusion product by the following lemma:

\begin{Lem}\label{Lem:Isomorphim_of_fusion}
  Assume that each generator $v_j$ of $M^j$ is annihilated by $\fn_0^-$.
  Then we have the following isomorphisms of $\fb$-modules:\\
  {\normalfont (i)} 
    \[ M^1_{0} * M^2_{c_2} *\cdots * M^p_{c_p}
       \cong \left[ M^1 * M^2_{c_2} * \cdots *M^p_{c_p}\right]_{\fb}, \ \text{and}
    \]
  {\normalfont(ii)} 
    \[ M^1_{c_1} * \cdots * M^p_{c_p} \cong \left[\C_{\mathrm{triv}} * M^1_{c_1} * \cdots * M^p_{c_p}\right]_\fb,
    \]
    where $\C_{\mathrm{triv}}$ denotes the trivial module.
\end{Lem}

\begin{proof}
  From the assumption, $M^1$ is generated by $v_1$ as a $\fb$-module.
  Hence the right hand side of the isomorphism (i) makes sense.
  The isomorphisms easily follow from the definition since we have
  \[ U(\fp_{I_0}')^{\le k}(v_1' \otimes v_2' \otimes \cdots \otimes v_p') = 
     U(\fb')^{\le k}(v_1' \otimes v_2' \otimes \cdots \otimes v_p')
  \]
  by the assumption and the Poincar\'{e}-Birkhoff-Witt theorem. 
\end{proof}

\begin{Rem}\normalfont
  The isomorphism (ii) of the above lemma does not hold in general.
  For example, let $M^1$ be a finite-dimensional irreducible $\fg_0$-module with a highest weight vector $v_1$,
  which is considered as a $\fp_{I_0}$-module via the evaluation map $\fp_{I_0} \to \fg_0\colon
  x \otimes t^k \mapsto \gd_{0,k}x, \
  K,d \mapsto 0$.
  Then by Lemma \ref{Lem:elementary} (iv), the fusion product $M_{c_1,v_1}^1$ is isomorphic to $M^1$.
  However, we easily see that the degree $0$ space of $\left[ \C_{\mathrm{triv}} * M^1_{c_1,v_1}\right]_{\fb}$ 
  is one-dimensional, and hence they are not isomorphic unless $M^1$ is trivial. 
  On the other hand, $\left[\C_{\mathrm{triv}} * M^1\right]_{\fb}$ is isomorphic to $M^1$
  if a lowest weight vector of $M^1$ is chosen as a generator.
  As seen from this example, the $\fb$-fusion product is sensitive to the choice of generators.
\end{Rem}

\begin{Lem}\label{Lem:important}
  Let $i \in I$. 
  The $\fb$-module $\left[ N * M^1 * \cdots * M^p\right]_{\fb}$ extends to a $\fp_i$-module if 
  $N$ extends to a $\fp_i$-module and either of the following conditions is satisfied:\\
  {\normalfont (i)} $i \in I_0$ and all $v_j$ and $u$ are annihilated by $f_i$, or \\
  {\normalfont (ii)} $i=0$, $K$ acts trivially on $M^1, \ldots ,M^p$, 
     and each $v_j$ {\normalfont(}resp.\ $u${\normalfont)} is annihilated by $e_\gt \otimes \C[t]$ 
     {\normalfont(}resp.\ $e_\gt \otimes t^{-1}${\normalfont)}.  
\end{Lem}

\begin{proof}
  The case (i) is easily proved since
  $N \otimes M_{c_1}^1 \otimes \cdots \otimes M_{c_p}^p$ is a $\fp_i'$-module and we have
  \begin{equation}\label{eq:coincidence}
     U(\fb')^{\le k}(u \otimes v_1' \otimes \cdots \otimes v_p') =
     U(\fp_i')^{\le k}(u \otimes v_1' \otimes \cdots \otimes v_p')
  \end{equation}
  for each $k$.
  Let us prove the case (ii).
  Since $K$ acts trivially, a $\fp'_0$-module structure is defined on each $M^j_{c_j}$ 
  by letting $f_0$ act by
  \[ \big(e_{\gt} \otimes t^{-1}\big)w' = - \sum_{k \ge 0} (-c_j)^{-k-1} \big((e_{\gt} \otimes t^k)w\big)' \ \ \ \text{for} \
     w\in M^j.
  \]
  Note that the above sum is finite since $M^j \in \Mod{\fp_{I_0}}$.
  Hence $N \otimes M_{c_1}^1 \otimes \cdots \otimes M_{c_p}^p$ extends to a $\fp_0'$-module.
  Moreover the equality (\ref{eq:coincidence}) also holds in this case from the assumption. Hence the assertion also follows
 in this case.  
\end{proof}

Take arbitrary vectors $w_j \in M^j$ ($1 \le j\le p$) and $z \in N$.
Let $k$ be the unique integer such that
\[ z \otimes w_1' \otimes \cdots \otimes w_p' \in
   \wti{F}_{u \otimes \cdots \otimes v_p'}^k (N\otimes M^1_{c_1} \otimes \cdots \otimes M^p_{c_p})
   \setminus \wti{F}_{u \otimes \cdots \otimes v_p'}^{k-1} (N\otimes M^1_{c_1} \otimes \cdots \otimes M^p_{c_p}),
\]
and denote by $z * w_1 * \cdots * w_p$ the vector of $\left[ N * M^1 * \cdots * M^p\right]_{\fb}$
which is the image of $z \otimes w_1' \otimes \cdots \otimes w_p'$ under the projection
\[ \wti{F}^{k}(N \otimes \cdots \otimes M^p_{c_p}) \twoheadrightarrow \wti{F}^{k}(N \otimes \cdots \otimes M^p_{c_p})
   / \wti{F}^{k-1}(N \otimes \cdots \otimes M^p_{c_p}).
\]
Note that $u * v_1 * \cdots *v_p$ is a generator of $\left[ N * M^1 * \cdots * M^p\right]_{\fb}$.
The following lemma, which obviously follows by definition, is important for the later arguments.

\begin{Lem}\label{Lem:annihilating}
  Let $X \in U(\fb')^{k}$. 
  Then $X$ annihilates $u * v_1 * \cdots * v_p$ if and only if there exists some $Y \in U(\fb')^{\le k-1}$ satisfying
  \[ (X-Y)(u \otimes v_1' \otimes \cdots \otimes v_p') = 0.
  \]
\end{Lem}

\section{Proof of the main theorem}\label{Section:proof}

Now, we begin the proof of Theorem \ref{Thm:Main_theorem}.
Assume that $\fg$ is a nontwisted affine Lie algebra and $\fg_0$ is of type $ADE$.
For a given sequence of KR modules $KR^{r_1,\ell_1},\ldots, KR^{r_p,\ell_p}$,
we set $M^j = KR^{r_j,\ell_j}$ and $v_j = v_{r_j,\ell_j} \in M^j$ for $1\le j \le p$ for short, 
and write $v_j^w= v_{r_j,\ell_j}^w \in M^j$ for $w \in W_0$ (defined in Section \ref{Section:KR}).

We shall show the theorem by the induction on $p$.
The assertion of the theorem for $p =1$ follows from Lemma \ref{Lem:elementary} (iv), Theorem \ref{Thm:construction},
Theorem \ref{Thm:KR-Demazure}, and Lemma \ref{Lem:trivial_lem2} (i).
Assume $p >1$. 
By Lemma \ref{Lem:elementary} (iii), we may (and do) assume $c_p = 0$, which implies $c_1, \ldots,c_{p-1}$ are nonzero.
First we show the following lemma:

\begin{Lem}\label{Lem:translation}
  We have the following isomorphisms of $\fb'$-modules:
  \begin{align}\label{eq:isom1}
     \C_{\ell_p\gL_0}\otimes \Big(M^p_{v_p} &* M^{p-1}_{v_{p-1}}* \cdots *M^1_{v_1}\Big) \nonumber\\
     &\cong\Big[ V_{t_{w_0(\varpi_{r_p})}}(\ell_p\gL_0)_{u_{t_{w_0(\varpi_{r_p})}(\ell_p\gL_0)}} 
      * M^{p-1}_{v_{p-1}^{w_0}} * \cdots *M^1_{v_1^{w_0}}\Big]_{\fb},
  \end{align}
  and 
  \begin{align}\label{eq:isom2}
    \cD_{t_{w_0(\varpi_{r_p})}}\Big(\C_{(\ell_p-\ell_{p-1})\gL_0} \otimes \, & \cD_{t_{w_0(\varpi_{r_{p-1}})}}
    \big(\C_{(\ell_{p-1}-\ell_{p-2}) \gL_0} \otimes \cdots \otimes \cD_{t_{w_0(\varpi_{r_1})}}\C_{\ell_1\gL_0}\big)\Big) 
    \nonumber \\
    &\cong \cD_{t_{w_0(\varpi_{r_p})}}\left[ \C_{\ell_p\gL_0} * M^{p-1}_{v_{p-1}^{w_0}} 
    * \cdots * M^1_{v_1^{w_0}}\right]_{\fb}.
  \end{align}
\end{Lem}

\begin{proof}
  Since
  \[ U(\fg_0)(v_p^{w_0} \otimes \cdots \otimes v_1^{w_0}) = U(\fg_0)(v_p \otimes \cdots \otimes v_1)
  \]
  holds, we have from Lemma \ref{Lem:elementary} (ii) that
  \begin{align*}
   M^{p}_{0,v_{p}} * M^{p-1}_{c_{p-1},v_{p-1}} * \cdots *M^1_{c_1,v_1}
   \cong M^{p}_{0,v_{p}^{w_0}} * M^{p-1}_{c_{p-1},v_{p-1}^{w_0}} * \cdots *M^1_{c_1,v_1^{w_0}},
  \end{align*}
  whose right hand side is isomorphic to 
  $\left[M^{p}_{v_{p}^{w_0}} * M^{p-1}_{v_{p-1}^{w_0}} *\cdots *M^1_{v_1^{w_0}}\right]_{\fb}$ by Lemma 
  \ref{Lem:Isomorphim_of_fusion} (i).
  Hence we have using (\ref{eq:level}) and Theorem \ref{Thm:KR-Demazure} that
  \begin{align*}
     \C_{\ell_p\gL_0}\otimes \Big(M^p_{v_p} * \cdots *&M^1_{v_1}\Big)
     \cong \left[ \big(\C_{\ell_p\gL_0} \otimes M^p\big)_{u_{\ell_p\gL_0} \otimes v_p^{w_0}} * M^{p-1}_{v_{p-1}^{w_0}}* 
     \cdots * M^1_{v_1^{w_0}} \right]_{\fb}\\
     & \ \cong\Big[ V_{t_{w_0(\varpi_{r_p})}}(\ell_p\gL_0)_{u_{t_{w_0(\varpi_{r_p})}(\ell_p\gL_0)}} 
      * M^{p-1}_{v_{p-1}^{w_0}} * \cdots *M^1_{v_1^{w_0}}\Big]_{\fb}.
  \end{align*}
  The isomorphism (\ref{eq:isom1}) is proved.
  Let us prove (\ref{eq:isom2}).
  By the induction hypothesis, there exists an isomorphism
  \begin{align*}
    \cD_{t_{w_0(\varpi_{r_{p-1}})}}\big(\C_{(\ell_{p-1} - \ell_{p-2})\gL_0} \otimes \cdots &
    \otimes \cD_{t_{w_0(\varpi_{r_1})}}\C_{\ell_1\gL_0}\big) \\
    &\cong \C_{\ell_{p-1} \gL_0} \otimes \big(M^{p-1}_{v_{p-1}} * \cdots *M^1_{v_1}\big),
  \end{align*}
  whose right hand side is isomorphic to 
  \begin{align*}
   \C_{\ell_{p-1} \gL_0} \otimes \big(M^{p-1}_{v_{p-1}^{w_0}} * \cdots *M^1_{v_1^{w_0}}\big)
   \cong \C_{\ell_{p-1} \gL_0} \otimes \Big[\C_{\mathrm{triv}}* M^{p-1}_{v_{p-1}^{w_0}} * \cdots *M^1_{v_1^{w_0}} \Big]_\fb
  \end{align*}
  by Lemma \ref{Lem:elementary} (ii) and Lemma \ref{Lem:Isomorphim_of_fusion} (ii).
  Hence we have using (\ref{eq:level}) that
  \begin{align*}
    \big(\text{left hand side} \  \text{of (\ref{eq:isom2}})\big)& \\ 
    \cong \cD_{t_{w_0(\varpi_{r_p})}}&\Big(\C_{\ell_p \gL_0}\otimes \left[\C_{\mathrm{triv}} * 
    M^{p-1}_{v_{p-1}^{w_0}} * \cdots *M^1_{v_1^{w_0}}\right]_{\fb}\Big) \\
    &\cong \cD_{t_{w_0(\varpi_{r_p})}}\left[\C_{\ell_p\gL_0} * 
    M^{p-1}_{v_{p-1}^{w_0}} * \cdots *M^1_{v_1^{w_0}}\right]_{\fb}.
  \end{align*}
  The isomorphism (\ref{eq:isom2}) is proved.
  \end{proof}
By Lemmas \ref{Lem:translation} and \ref{Lem:trivial_lem2} (i), in order to prove the theorem
it suffices to show the following isomorphism of $\fb'$-modules:
\begin{align}\label{eq:translated_isom}
  \Big[ V_{t_{w_0(\varpi_{r_p})}}(\ell_p\gL_0&)_{u_{t_{w_0(\varpi_{r_p})}(\ell_p\gL_0)}}
      * M^{p-1}_{v_{p-1}^{w_0}} * \cdots *M^1_{v_1^{w_0}}\Big]_{\fb} \nonumber \\
  & \ \ \ \ \ \cong \cD_{t_{w_0(\varpi_{r_p})}}\Big[ \C_{\ell_p\gL_0} * M^{p-1}_{v_{p-1}^{w_0}} * 
  \cdots * M^1_{v_1^{w_0}}\Big]_{\fb}.
\end{align}
Let $w \in W$ and $\tau \in \gS$ be unique elements satisfying $w\tau = t_{w_0(\varpi_{r_p})}$, and 
$w=s_{i_k} \cdots s_{i_1}$ a reduced expression. 
For $0 \le q \le k$, let $w^q = s_{i_q} \cdots s_{i_1}\tau \in \wti{W}$, and $\ol{w}^q$ be the unique element of $W_0$ 
satisfying $w^q \in \ol{w}^q\cdot T(\wti{M})$ (note that $\ol{w}^k = \id$).
We also write $\ol{\tau}$ for $\ol{w}^0$.
Then since 
\[ V_{t_{w_0(\varpi_{r_p})}}(\ell_p\gL_0) \cong \cD_{t_{w_0(\varpi_{r_p})}}\C_{\ell_p\gL_0} 
   \cong \cD_{w}\C_{\ell_p\gL_{\tau(0)}}
\]
holds, we see that the isomorphism (\ref{eq:translated_isom}) is deduced 
by the induction on $q$ from the following two propositions, and hence Theorem \ref{Thm:Main_theorem} is established:

\begin{Prop}\label{Prop1}
  We have
  \begin{align*}
    \Big[ \C_{\ell_p\gL_{\tau(0)}} * M_{c_{p-1},v_{p-1}^{\ol{\tau}w_0}}^{p-1} &
    * \cdots * M_{c_1,v_1^{\ol{\tau}w_0}}^{1}\Big]_{\fb} \\
    &\cong (\tau^{-1})^*\left[\C_{\ell_p\gL_0} * M_{c_{p-1},v_{p-1}^{w_0}}^{p-1} * \cdots *M_{c_1,v_1^{w_0}}^1 \right]_{\fb}
  \end{align*}    
  as $\fb'$-modules.
\end{Prop}

\begin{Prop}\label{Prop2}
  For each $1 \le q \le k$, we have
  \begin{align*}
    \Big[ V_{w^{q}}(\ell_p\gL_0&)_{u_{w^q(\ell_p\gL_0)}}
    * M_{v_{p-1}^{\ol{w}^qw_0}}^{p-1} * \cdots * M_{v_1^{\ol{w}^qw_0}}^{1}\Big]_{\fb}\\
    &\cong \cD_{i_q}\Big[ V_{w^{q-1}}(\ell_p\gL_0)_{u_{w^{q-1}(\ell_p\gL_0)}}
    * M_{v_{p-1}^{\ol{w}^{q-1}w_0}}^{p-1} * \cdots * M_{v_1^{\ol{w}^{q-1}w_0}}^1 \Big]_{\fb}
  \end{align*}
  as $\fb'$-modules.
\end{Prop}

To show Proposition \ref{Prop1}, we need to prepare several lemmas.
The following one is proved similarly as \cite[Lemma 3.8]{MR1104219}:

\begin{Lem}\label{Lem:Kac}
  Let $M$ be a finite-dimensional $\fp_{I_0}'$-module.
  Then for every $w \in W_0$, there exists a linear automorphism $\eta_w$ of $M$ satisfying
  \begin{align*}
    \Ad(\eta_w) (h \otimes t^s)= &w(h) \otimes t^s \ \ \text{for} \ h \in \fh_0, s \in \Z_{\ge 0}, \ \ \ \Ad(\eta_w)(K) =K, \\
    \Ad(\eta_w) (e_\ga \otimes t^s) &= a_{w}(\ga)e_{w(\ga)} \otimes t^s \ \ \text{for}\
     \ga \in \gD_0, s \in \Z_{\ge 0},
     \ \ \ \text{and} \\
    &\eta_w(M_{\gl}) =M_{w(\gl)} \ \text{for} \ \gl \in \fh^*/\C\gd,
  \end{align*}
  where $a_w(\ga)$ are some nonzero complex numbers which do not depend on $M$.
\end{Lem}

By applying $\Ad(\eta_{w_0})$ given in Lemma \ref{Lem:Kac} to the defining relations of $KR^{r,\ell}$ in Definition
\ref{Def:KR_module}, the following lemma is proved:

\begin{Lem}\label{Lem:1-0}
  The annihilating ideal of $v_{r,\ell}^{w_0} \in KR^{r,\ell}$ in $U\big(\fp_{I_0}'\big)$ is generated by 
  \begin{align*}
    \fn^-_0 \otimes\, & \C[t], \ \ \ \fh_0 \otimes t\C[t],\ \ \ h- \langle w_0(\ell\varpi_r) ,h \rangle \ \ (h \in \fh'), \\
    & e_{\bar{r}}^{\ell+1},\ \ \ e_{\bar{r}}\otimes t, \ \ \ \text{and} \ \ e_i \ \ (i \in I_0\setminus\{\bar{r}\}),
  \end{align*}
  where $\bar{r}$ is the node of $I_0$ such that $w_0(\ga_r) = - \ga_{\bar{r}}$.
\end{Lem}

\begin{Lem}\label{Lem:1-1}
  For $c \in \C$, the annihilating ideal of ${v_{r,\ell}^{w_0}}' \in KR^{r,\ell}_c$ in $U(\fb')$ is generated by
  \begin{align*} 
     \fn^-_0 \otimes \, & t\C[t],\ \ \ \fh_0 \otimes (t-c)\C[t],  
     \ \ \ h- \langle w_0(\ell\varpi_r) ,h \rangle \ \ (h \in \fh'),\\
     &e_{\bar{r}}^{\ell+1}, \ e_{\bar{r}} \otimes (t-c), \ \ \ \text{and}\ \ e_i \ \ (i \in I_0 \setminus \{\bar{r}\}).
  \end{align*}
\end{Lem}

\begin{proof}
  Let $I$ be the subspace of $U(\fb')$ spanned by the above vectors.
  From Lemma \ref{Lem:1-0},
  we see that the annihilating ideal of ${v_{r,\ell}^{w_0}}'$ in $U\big(\fp_{I_0}'\big)$ is equal to  
  $U\big(\fp_{I_0}'\big)\big(I + \fn^-_0 \big)$.
  We have to prove that
  \[ U\big(\fp_{I_0}'\big)\big(I + \fn^-_0 \big) \cap U(\fb') \subseteq U(\fb')I.
  \]
  Since $U\big(\fp_{I_0}'\big) = U(\fb') \oplus U\big(\fp_{I_0}'\big)\fn_0^-$ holds
  by the Poincar\'{e}-Birkhoff-Witt theorem, it suffices to show that
  \[ U(\fp_{I_0}')I \subseteq U(\fb')I \oplus U\big(\fp_{I_0}'\big)\fn^-_0.
  \]
  Then since we have $U(\fp_{I_0}') = U(\fb')U(\fn^-_0)$ and 
  $U(\fn^-_0)$ is generated by $\{f_i \mid i \in I_0\}$, it is enough to prove that 
  $f_i I \subseteq I \oplus U\big(\fp_{I_0}'\big)\fn^-_0$ holds for $i \in I_0$,  
  which is proved by elementary calculations. 
\end{proof}

\begin{Lem}\label{Lem:annihilators}
  There exists a nonzero complex number $b$ satisfying the following statement:
  for every KR module $KR^{r,\ell}$ and $c \in \C$, 
  there exists a homomorphism 
  \[ KR^{r,\ell}_c \to \tau^* KR^{r,\ell}_{bc}
  \]
  of $\fb'$-modules which maps ${v_{r,\ell}^{w_0}}'$ to $\tau^*{v_{r,\ell}^{\ol{\tau}w_0}}'$.
\end{Lem}

\begin{proof}
  Since $KR^{r,\ell}_c$ is generated by ${v_{r,\ell}^{w_0}}'$ as a $\fb'$-module,
  it suffices to show for suitable $b \in \C^*$ that if $X \in U(\fb')$ annihilates ${v_{r,\ell}^{w_0}}' \in KR^{r,\ell}_{c}$,
  then $\tau(X)$ annihilates ${v_{r,\ell}^{\ol{\tau}w_0}}' \in KR^{r,\ell}_{bc}$.
  Since the $\fh'$-weights of ${v_{r,\ell}^{w_0}}'\in KR^{r,\ell}_c$ 
  and $\tau^*{v_{r,\ell}^{\ol{\tau}w_0}}'\in \tau^*KR^{r,\ell}_{bc}$ coincide,
  we may assume $X \in U(\fn^+)$.
  Let $i_0 = \tau^{-1}(0) \in I$, and $\psi_{i_0}\colon U(\fb') \to U\big(\fp_{I_0}'\big)$ be an algebra homomorphism defined by
  \begin{align*}
    &\psi_{i_0}(x \otimes t^s) = x \otimes t^{s + \langle \ga ,\varpi_{i_0}^\vee \rangle } \ \ \ \text{for} \ \ x \in \fg_{\ga}
    \ (\ga \in \gD_0),\\
    & \psi_{i_0}(h \otimes t^s) = h \otimes t^s \ \ \ \text{for} \ \ h \in \fh_0, \ \ \ \psi_{i_0}(K) = K.
  \end{align*}
  Using Lemma \ref{Lem:1-1}, we easily check that $\psi_{i_0}(X)$ also annihilates ${v_{r,\ell}^{w_0}}'\in KR^{r,\ell}_c$.
  Let $\eta= \eta_{\ol{\tau}}$ be the linear automorphism of $KR^{r,\ell}_c$ given in Lemma \ref{Lem:Kac}.
  Then since $\eta({v_{r,\ell}^{w_0}}') \in \C^*{v_{r,\ell}^{\ol{\tau}w_0}}'$ holds by Lemma \ref{Lem:KR},
  $\Ad(\eta) \circ \psi_{i_0}(X)$ annihilates ${v_{r,\ell}^{\ol{\tau}w_0}}' \in KR^{r,\ell}_c$, 
  which is equivalent to that $\varphi_c \circ \Ad(\eta) \circ \psi_{i_0}(X)$ annihilates 
  $v_{r,\ell}^{\ol{\tau}w_0} \in KR^{r,\ell}$ ($\varphi_c$ is defined in Section \ref{Fusion product}).
  It is easy to see from Lemma \ref{Lem:transform} that there exists some $b_i \in \C^*$ for each $i \in I$ such that
  \[ \Ad(\eta) \circ \psi_{i_0}(e_{\tau^{-1}(i)}) =  b_i e_i.
  \]
  Define a linear automorphism $H$ on $KR^{r,\ell}$ by
  \[ H(u) = \prod_{i \in I} b_i^{- \langle \gl, \gL_i^{\vee}\rangle}\cdot u \ \ \ 
     \text{if} \ u \in KR^{r,\ell}_{\gl} \ \ (\gl \in P),
  \]
  where $\gL_i^{\vee}\in \fh$ are the fundamental coweights of $\fg$. 
  $\Ad(H)\circ\varphi_c\circ\Ad(\eta)\circ\psi_{i_0}(X)$ annihilates $v_{r,\ell}^{\ol{\tau}w_0}$ since it is a weight vector.
  Set $b = \prod_{i \in I} b_i^{a_i}$. 
  It is easily checked that
  \[ \Ad(H)\circ \varphi_c\circ\Ad(\eta)\circ\psi_{i_0}(e_{\tau^{-1}(i)}) = \begin{cases} e_{-\gt} \otimes 
                                                                                                          t + bc e_{-\gt} & 
                                                                                     \text{if} \ i = 0, \\
                                                                                     e_{i} & \text{otherwise},
                                                                       \end{cases}
  \]
  which implies
  \[ \Ad(H)\circ \varphi_c\circ\Ad(\eta)\circ\psi_{i_0} = \varphi_{bc} \circ \tau \ \ \ \text{on} \ \ U(\fn^+).
  \] 
  Hence we see that $\varphi_{bc} \circ \tau(X)$ annihilates $v_{r,\ell}^{\ol{\tau}w_0}$, 
  which is equivalent to that
  $\tau(X)$ annihilates ${v_{r,\ell}^{\ol{\tau}w_0}}' \in KR^{r,\ell}_{bc}$. The assertion is proved.
\end{proof}

Now, we give the proof of Proposition \ref{Prop1}:\\

\noindent \textit{Proof of Proposition \ref{Prop1}.\ }
  Note that the right hand side of Theorem \ref{Thm:Main_theorem} does not depend on the parameters $c_1,\ldots,c_p$.
  Hence from the induction hypothesis on $p$ and the proof of (\ref{eq:isom2}), we see that the $\fb$-module 
  $\Big[\C_{\ell_p\gL_0} * M_{c_{p-1},v_{p-1}^{w_0}}^{p-1} * \cdots * M_{c_1,v_1^{w_0}}^1\Big]_{\fb}$
  also does not depend on the parameters, and in particular we have
  \[ \Big[\C_{\ell_p\gL_0} * M_{c_{p-1},v_{p-1}^{w_0}}^{p-1} * \cdots * M_{c_1,v_1^{w_0}}^1\Big]_{\fb} \cong
     \Big[\C_{\ell_p\gL_0} * M_{b^{-1}c_{p-1},v_{p-1}^{w_0}}^{p-1} * \cdots * M_{b^{-1}c_1,v_1^{w_0}}^1\Big]_{\fb},
  \]
  where $b$ is the complex number given in Lemma \ref{Lem:annihilators}.
  Hence the proposition is equivalent to the following isomorphism of $\fb'$-modules:
  \begin{align*}
    \Big[ \C_{\ell_p\gL_{\tau(0)}} * M_{c_{p-1},v_{p-1}^{\ol{\tau}w_0}}^{p-1} &
    * \cdots * M_{c_1,v_1^{\ol{\tau}w_0}}^{1}\Big]_{\fb} \\
    &\cong (\tau^{-1})^*\left[\C_{\ell_p\gL_0} * M_{b^{-1}c_{p-1},v_{p-1}^{w_0}}^{p-1} * \cdots 
    *M_{b^{-1}c_1,v_1^{w_0}}^1 \right]_{\fb}.
  \end{align*}    
  Let us prove this.
  Let $u_1 \in \C_{\ell_p \gL_0}$ and $u_2 \in \C_{\ell_p \gL_{\tau(0)}}$ be nonzero vectors.
  Since dimensions of two modules are equal, it suffices to show there exists a surjective 
  homomorphism of $\fb'$-modules from the right hand side to the left hand side
  mapping $(\tau^{-1})^*(u_1 * v_{p-1}^{w_0} * \cdots * v_1^{w_0})$ to 
  $u_2 * v_{p-1}^{\ol{\tau}w_0} *\cdots * v_1^{\ol{\tau} w_0}$,
  which is equivalent to show that if $X \in U(\fb')$ annihilates $u_1 * v_{p-1}^{w_0} * \cdots * v_1^{w_0}$, 
  then $\tau(X)$ annihilates $u_2 * v_{p-1}^{\ol{\tau}w_0} *\cdots * v_1^{\ol{\tau} w_0}$.
  We may assume $X \in U(\fb')_{\gg}^s$ for some $\gg \in Q_0$ and $s \in \Z_{\ge 0}$, where we set
  \[ U(\fb')_{\gg}^s = \left\{ Z \in U(\fb')^s \mid [h,Z] = \langle \gg, h \rangle Z \ \text{for} \ h \in \fh_0 \right\}.
  \]  
  By Lemma \ref{Lem:annihilating}, there exists $Y \in U(\fb')^{\le s-1}$ such that
  \begin{equation}\label{eq:trivial}
    (X-Y)(u_1 \otimes {v_{p-1}^{w_0}}' \otimes \cdots \otimes {v_1^{w_0}}') 
    = 0,
  \end{equation}
  and we may assume $Y\in U(\fb')^{\le s-1}_\gg$.
  We see from Lemma \ref{Lem:annihilators} that there exists a homomorphism of $\fb'$-modules 
  \[ \C_{\ell_p\gL_0} \otimes M^{p-1}_{b^{-1}c_{p-1}} \otimes \cdots \otimes M^1_{b^{-1}c_1} \to 
     \tau^*\Big(\C_{\ell_p \gL_{\tau(0)}} \otimes M^{p-1}_{c_{p-1}} \otimes \cdots \otimes M^1_{c_1}\Big)
  \]
  which maps $u_1 \otimes {v_{p-1}^{w_0}}' \otimes \cdots \otimes {v_1^{w_0}}'$
  to $\tau^*\big(u_2 \otimes {v_{p-1}^{\ol{\tau}w_0}}' \otimes \cdots \otimes {v_1^{\ol{\tau} w_0}}'\big)$.
  From this and (\ref{eq:trivial}), we have
  \[ \tau(X - Y)(u_2 \otimes {v_{p-1}^{\ol{\tau}w_0}}' \otimes \cdots \otimes {v_1^{\ol{\tau} w_0}}') = 0.
  \]
  Moreover, we have $\tau(X) \in U(\fb')_{\ol{\tau}(\gg)}^{s + \langle \gg, \varpi_{i_0}^\vee\rangle}$ and 
  $\tau(Y) \in U(\fb')_{\ol{\tau}(\gg)}^{\le s +\langle \gg, \varpi_{i_0}^\vee\rangle-1}$ from Lemma \ref{Lem:transform}, 
  where we set $i_0 = \tau^{-1}(0)$.
  Hence by Lemma \ref{Lem:annihilating}, 
  $\tau(X)$ annihilates $u_2 * v_{p-1}^{\ol{\tau}w_0} * \cdots * v_1^{\ol{\tau}w_0}$. 
  The assertion is proved. \qed \\

It remains to prove Proposition \ref{Prop2}:\\

\noindent \textit{Proof of Proposition \ref{Prop2}.\ }
  We abbreviate $u^q = u_{w^q(\ell_p\gL_0)}$, $V_q = V_{w^q}(\ell_p \gL_0)$ and $v_{j}^{q} = v_{j}^{\ol{w}^qw_0}$.
  Fix $1\le q\le k$, and assume when $q > 1$ that the assertion of the proposition holds for $q' < q$.
  Then we see that $\Big[ (V_{q-1})_{u^{q-1}} * M_{v_{p-1}^{q-1}}^{p-1} * \cdots * M_{v_1^{q-1}}^1\Big]_{\fb}$
  has a Demazure flag from the induction hypothesis on $p$, Corollary \ref{Cor:Joseph}, Lemma \ref{Lem:affine_version} (ii),
  and Proposition \ref{Prop1}.

  Recall that, by the definition of Demazure modules, there exists a canonical embedding $V_{q-1} \hookrightarrow V_q$.
  Let $z^{q-1}$ denote the image of $u^{q-1}$ under this embedding. 
  It should be noted that $z^{q-1} * v_{p-1}^{q-1} * \cdots * v_1^{q-1} 
  \in \Big[ (V_q)_{u^q} * M_{v_{p-1}^{q}}^{p-1} * \cdots * M_{v_1^q}^1\Big]_{\fb}$ is not a generator. \\ 
  \\ \textit{Claim} 1.\ \ There exists a homomorphism of $\fb'$-modules 
    \[ \Phi\colon\Big[ (V_{q-1})_{u^{q-1}} * M_{v_{p-1}^{q-1}}^{p-1} * \cdots * M_{v_1^{q-1}}^1\Big]_{\fb} \to
       \Big[ (V_q)_{u^q} * M_{v_{p-1}^{q}}^{p-1} * \cdots * M_{v_1^q}^1\Big]_{\fb}
    \]
    which maps $u^{q-1} * v_{p-1}^{q-1} * \cdots * v_1^{q-1}$ to $z^{q-1} * v_{p-1}^{q-1} * \cdots* v_1^{q-1}$. \\[-5pt]
    
    It suffices to show that, if $X \in U(\fb')$ annihilates $u^{q-1} * v_{p-1}^{q-1} * \cdots * v_1^{q-1}$,
    then $X$ also annihilates $z^{q-1} * v_{p-1}^{q-1} * \cdots * v_1^{q-1}$.
    We may assume that $X \in U(\fb')^s$ for some $s \in \Z_{\ge 0}$.
    Then there exists some $Y \in U(\fb')^{\le s-1}$ satisfying 
    \[ (X-Y)(u^{q-1} \otimes {v_{p-1}^{q-1}}' \otimes \cdots \otimes {v_1^{q-1}}') =0
    \]
    by Lemma \ref{Lem:annihilating}.
    Obviously,
    \begin{align}\label{eq:X-Y}
      (X-Y)(z^{q-1} \otimes {v_{p-1}^{q-1}}' \otimes \cdots \otimes {v_1^{q-1}}' )= 0
    \end{align}
    also holds.
    Let $N$ be the unique integer such that
    \begin{align*}
      z^{q-1} \otimes {v_{p-1}^{q-1}}' \otimes \cdots \otimes {v_1^{q-1}}' 
      &\notin U(\fb')^{\le N-1} (u^q \otimes {v_{p-1}^q}' \otimes \cdots \otimes {v_1^q}') \ \text{and} \\
      z^{q-1} \otimes {v_{p-1}^{q-1}}' \otimes \cdots \otimes {v_1^{q-1}}' &\in 
      U(\fb')^{\le N}(u^q \otimes {v_{p-1}^q}' \otimes \cdots \otimes {v_1^q}'),
    \end{align*}
    and $Z_N \in U(\fb')^{N}$ and $Z_{\le N-1} \in U(\fb')^{\le N-1}$ be vectors such that 
    \[ (Z_N + Z_{\le N-1})(u^q \otimes {v_{p-1}^q}' \otimes \cdots \otimes {v_1^q}') 
       = z^{q-1} \otimes {v_{p-1}^{q-1}}' \otimes \cdots \otimes {v_1^{q-1}}'.
    \]
    Then we have from (\ref{eq:X-Y}) that
    \begin{align*}
      XZ_N(&u^q \otimes {v_{p-1}^q}' \otimes \cdots \otimes {v_1^q}')\\
      &= X(z^{q-1} \otimes {v_{p-1}^{q-1}}' \otimes \cdots \otimes {v_1^{q-1}}') 
         -XZ_{\le N-1}(u^q \otimes {v_{p-1}^q}' \otimes \cdots \otimes {v_1^q}')\\
       &= \big(Y(Z_N + Z_{\le N-1}) - XZ_{\le N-1}\big)(u^q \otimes {v_{p-1}^q}' \otimes \cdots \otimes {v_1^q}').
    \end{align*}
    Since 
    \[ XZ_N \in U(\fb')^{s +N} \ \ \text{and} \ \ Y(Z_N + Z_{\le N-1}) - XZ_{\le N-1} \in U(\fb')^{\le s+N-1},
    \]
    $XZ_N(u^q * v_{p-1}^q * \cdots * v_1^q) = 0$ holds by Lemma \ref{Lem:annihilating}.
    On the other hand, we have by definition that
    \[ Z_N(u^q* v_{p-1}^q * \cdots * v_1^q) = z^{q-1} * v_{p-1}^{q-1} * \cdots * v_1^{q-1}.
    \]    
    Hence
    \[ X(z^{q-1} * v_{p-1}^{q-1} * \cdots * v_1^{q-1}) = 0
    \]
    holds, and Claim 1 is proved. \\
    
  Set $i = i_q$. By Lemma \ref{Lem:important}, the $\fb$-module 
  $\Big[ (V_q)_{u^q} * M_{v_{p-1}^{q}}^{p-1} * \cdots * M_{v_1^q}^1\Big]_{\fb}$ extends to a $\fp_i$-module.
  Then by Lemma \ref{Lem:trivial_lem2} (ii),
  there exists a homomorphism of $\fp_i'$-modules
  \[ \wti{\Phi}\colon \cD_i\Big[ (V_{q-1})_{u^{q-1}} * M_{v_{p-1}^{q-1}}^{p-1} * 
   \cdots * M_{v_1^{q-1}}^1\Big]_{\fb} \to \Big[ (V_q)_{u^q} * M_{v_{p-1}^{q}}^{p-1} * 
   \cdots * M_{v_1^q}^1\Big]_{\fb}
  \]
  which makes the following diagram commutative:
  \[ \xymatrix{ \left[ (V_{q-1})_{u^{q-1}} * \cdots * M_{v_1^{q-1}}^1\right]_{\fb} \ar[r]^(.52){\Phi}
     \ar[d] & \left[ (V_q)_{u^q} * \cdots * M_{v_1^q}^1\right]_{\fb}  \\
     \cD_i\left[ (V_{q-1})_{u^{q-1}} *\cdots * M_{v_1^{q-1}}^1\right]_{\fb} \ar[ru]^(.50){\wti{\Phi}} %& 
     },
  \]  
  where the vertical map is the canonical one.\\
  \\ \textit{Claim} 2. \ \ The homomorphism $\wti{\Phi}$ is surjective. \\[-5pt]
  
    It suffices to show the image of $\wti{\Phi}$ contains the generator $u^q * v_{p-1}^q * \cdots * v_1^q$,
    whose $\fh'$-weight $\gl$ is equal to 
    \[ \gl = \cl\Big(w^q (\ell_p\gL_0) + \sum_{1 \le j \le p-1} \ol{w}^q w_0(\ell_j\varpi_{r_j})\Big) \in P_{\cl}.
    \] 
    Note that the image of $\wti{\Phi}$ contains $z^{p-1} * v_{p-1}^{q-1} * \cdots * v_1^{q-1}$, whose $\fh'$-weight $\mu$ is
    \[ \mu = \cl\Big(w^{q-1}(\ell_p \gL_0) + \sum_{1 \le j \le p-1} \ol{w}^{q-1} w_0(\ell_j\varpi_{r_j})\Big) = s_i(\gl).
    \]
    As the image is a $\fp_i'$-module, its weight set contains $s_i(\mu) = \gl$.
    Since $\Big[ (V_q)_{u^q} * \cdots * M_{v_1^q}^1\Big]_{\fb}$ is isomorphic to $V_q \otimes  \cdots 
    \otimes M^1$ as a $(\fg_0\oplus \C K)$-module,
    it is easily checked that the weight space with weight $\gl$ is one-dimensional.
    Hence the image contains the generator, and Claim 2 is proved.\\
  
  Now, the following claim completes the proof of the proposition:\\
  \\ \textit{Claim} 3. \ \ The dimensions of the both sides of $\wti{\Phi}$ are equal. \\[-5pt]
  
    The dimension of the right hand side is equal to
    \begin{equation}\label{eq:dimension}
      \dim V_q \times \prod_{1 \le j \le p-1} \dim M^j.
    \end{equation}
    Let us calculate the dimension of the left hand side. 
    As stated at the beginning of this proof, 
    $\Big[ V_{q-1} * M^{p-1} * \cdots * M^1\Big]_{\fb}$
    has a Demazure flag. 
    Hence by Corollary \ref{Cor:flag_character}, the character of the left hand side of $\wti{\Phi}$ is equal to
    \[ D_i\,  \ch \Big[ V_{q-1} * M^{p-1} * \cdots * M^1\Big]_{\fb}.
    \]
    For a $\fh'$-semisimple module $M$ whose $\fh'$-weight set is contained in $P_\cl$, 
    denote by $\bch M \in \Z[P_\cl]$ the $\fh'$-character of $M$.
    Since each $M^j$ is a finite-dimensional ($\fg_0\oplus \C K$)-module on which $K$ acts trivially, $\bch M^j$
    belongs to $\Z[P_{\cl}^0]$ and is $W_0$-invariant. 
    Hence we have from Lemma \ref{Lem:final_lem} that
    \begin{align*}
      \cl \circ D_i\,\ch \Big[ V_{q-1} * M^{p-1} * \cdots * M^1\Big]_{\fb}
      & = \ol{D}_i \, \bch \Big[ V_{q-1} * M^{p-1} * \cdots * M^1\Big]_{\fb} \\
        = \ol{D}_i \, \bigg(\bch V_{q-1} \times \prod_{1\le j \le p-1} \bch M^j\bigg)
      & = \prod_{1\le j \le p-1} \bch{M^j} \times \ol{D}_i \,\bch V_{q-1} \\
      & = \prod_{1 \le j \le p-1} \bch M^j \times \cl \circ D_i \, \ch V_{q-1}.
    \end{align*}
    Since $D_i\, \ch V_{q-1} = \ch V_q$ by Theorem \ref{Thm:Character_formula},
    we see that the dimension of the left hand side is equal to (\ref{eq:dimension}).
    Hence Claim 3 is proved, and the proof of the proposition is complete. \qed \\

As stated above, Theorem \ref{Thm:Main_theorem} is now established from Propositions \ref{Prop1} and \ref{Prop2}.

\begin{Rem} \normalfont
  In Theorem \ref{Thm:Main_theorem}, it is assumed that $\fg_0$ is of $ADE$ type.
  The author, however, expects the theorem to be true for general types if all given KR modules 
  satisfy the assumption of Theorem \ref{Thm:KR-Demazure}.
  In fact, all the proof of the theorem can also be applied in this case,
  except for the Joseph's theorem (Theorem \ref{Thm:Joseph}) which is needed in the final step of the proof.
  Hence to prove the theorem for non-simply laced type by our approach, it is needed to prove the Joseph's theorem for this type.
  Since Proposition \ref{Prop:X-side} has already been proved for nonexceptional type,
  this would also imply the $X=M$ conjecture for perfect KR crystals of type $B_n^{(1)}$ and $C_n^{(1)}$.  
\end{Rem}

%\bibliographystyle{plain}
%\bibliography{bibliography-Fermionic_form}

\def\cprime{$'$}

\end{document}